\documentclass[journal]{IEEEtran}
 
\IEEEoverridecommandlockouts                              

\usepackage{amsmath}    
\usepackage{amsfonts}
\usepackage{graphicx}   
\usepackage{subcaption}
\usepackage{epsfig} 
\usepackage{ragged2e}
\usepackage{color}
\usepackage[normalem]{ulem}
\usepackage{cancel}
\usepackage{amssymb}
\usepackage{color}
\usepackage{float}
\usepackage{cite}
\usepackage{enumitem}
\setlist[itemize]{leftmargin=*}

\newcommand{\MCnote}{\textcolor{black}}
\newcommand{\SHnote}{\textcolor{black}}

\newcommand{\R}{\mathbb{R}} 
\newcommand{\state}{z} 
\newcommand{\scstate}{x} 
\newcommand{\scctrl}{w} 
\newcommand{\spart}{\state} 
\newcommand{\cpart}{\ctrl} 

\newcommand{\traj}{\zeta} 

\newcommand{\ctrl}{u} 
\newcommand{\cset}{\mathcal{U}}

\newcommand{\cfset}{\mathbb{U}}

\newcommand{\dstb}{d} 
\newcommand{\dpart}{\dstb} 
\newcommand{\dset}{\mathcal{D}}
\newcommand{\scdstb}{b}
\newcommand{\dfset}{\mathbb{D}}

\newcommand{\genset}{\mathcal{S}}
\newcommand{\xset}{\mathcal{X}}
\newcommand{\zset}{\mathcal{Z}}

\newcommand{\proj}{\text{proj}}
\newcommand{\bp}{\proj^{-1}}

\newcommand{\fdyn}{f} 
\newcommand{\scdyn}{g} 
\newcommand{\sctraj}{\xi}

\newcommand{\fc}{l} 

\newcommand{\minbrs}{\mathcal A} 
\newcommand{\maxbrs}{\mathcal R} 
\newcommand{\minbrt}{\bar\minbrs} 
\newcommand{\maxbrt}{\bar\maxbrs} 
\newcommand{\targetset}{\mathcal T} 
\newcommand{\sctarget}[1]{\targetset_{#1}}
\newcommand{\pos}{p} 
\newcommand{\vel}{v} 

\newtheorem{rem}{Remark}
\newtheorem{defn}{Definition}
\newtheorem{prop}{Proposition}
\newtheorem{thm}{Theorem}
\newtheorem{lem}{Lemma}

\newtheorem{cor}{Corollary}

\title{\large \bf Decomposition of Reachable Sets and Tubes for a Class of Nonlinear Systems}

\author{Mo Chen, Sylvia L. Herbert, Mahesh S. Vashishtha, Somil Bansal, Claire J. Tomlin
\thanks{This work has been supported in part by NSF under CPS:ActionWebs (CNS-931843), by ONR under the HUNT (N0014-08-0696) and SMARTS (N00014-09-1-1051) MURIs and by grant N00014-12-1-0609, by AFOSR under the CHASE MURI (FA9550-10-1-0567). The research of M. Chen has received funding from the ``NSERC PGS-D'' Program.}
\thanks{All authors are with the Department of Electrical Engineering and Computer Sciences, University of California, Berkeley. \{mochen72, sylvia.herbert, msv, somil, tomlin\}@berkeley.edu}
\thanks{Theorem \ref{thm:sc_reach_in} of this paper is partially taken from our conference paper \cite{Chen17b}.}
}

\begin{document}
\maketitle
\thispagestyle{empty}
\pagestyle{empty}

\begin{abstract} 
Reachability analysis provides formal guarantees for performance and safety properties of nonlinear control systems. Here, one aims to compute the backward reachable set (BRS) or tube (BRT) -- the set of states from which the system can be driven into a target set at a particular time or within a time interval, respectively. The computational complexity of current approaches scales exponentially, making application to high-dimensional systems intractable. We propose a technique that decomposes the dynamics of a general class of nonlinear systems into subsystems which may be coupled through common states, controls, and disturbances. Despite this coupling, BRSs and BRTs can be computed efficiently using our technique without incurring additional approximation errors and without the need for linearizing dynamics or approximating sets as polytopes. Computations of BRSs and BRTs now become orders of magnitude faster, and for the first time BRSs and BRTs for many high-dimensional nonlinear control systems can be computed using the Hamilton-Jacobi (HJ) formulation. In situations involving bounded adversarial disturbances, our proposed method can obtain slightly conservative results. We demonstrate our theory by numerically computing BRSs and BRTs using the HJ formulation for several systems, including the 6D Acrobatic Quadrotor and the 10D Near-Hover Quadrotor.
\end{abstract}

\section{Introduction}
Many important real-world systems are described by complex nonlinear models whose behavior can be non-intuitive and difficult to predict. These systems include power \cite{Zhao2014,Dobbe2016}, biological \cite{Ghosh2001, Ghosh}, and robotic systems such as autonomous cars and unmanned aerial vehicles \cite{Hoffmann2007, Kong2015}. With the recent advancements of sophisticated system modeling in these areas, the dimensionality of system models has grown significantly. Many of these systems are also safety-critical, making their verification extremely important. As a result, computationally tractable tools for the analysis of these nonlinear, high-dimensional, and safety-critical systems are urgently needed.

Verification of systems is challenging for many reasons. First, all possible system behaviors must be accounted for. This makes most simulation-based approaches insufficient, and thus formal verification methods are needed. Second, many practical systems are affected by disturbances in the environment. In addition, the systems evolve in continuous time with complex, nonlinear dynamics. Lastly, perhaps the most difficult challenge of all is that these systems often have high-dimensional state spaces.

\MCnote{Reachability analysis is an important formal verification method for guaranteeing performance and safety properties of systems. However, current methods do not simultaneously address all of the above challenges. For example, \cite{Kong2015a, Duggirala2015} excel in determining whether system trajectories from a small set of initial conditions could potentially enter a set of unsafe states, but do not provide the backward reachable set (BRS) or tube (BRT) -- the set of all initial states from which entering some target set is inevitable. Due to the challenges of computing BRSs and BRTs, the state-of-the-art methods need to make trade offs on different axes of considerations such as computational scalability, generality of system dynamics, existence of control and/or disturbance variables, and flexibility in representation of sets.}

\MCnote{For example, the methods presented in \cite{Greenstreet1998, Frehse2011, Kurzhanski00, Kurzhanski02, Maidens13} have had success in analyzing relatively high-dimensional affine systems using sets of pre-specified shapes, such as polytopes or hyperplanes. Other potentially less scalable methods are able to handle systems with the more complex dynamics \cite{Chen2013, Althoff2015, Frehse2011, Majumdar13, Dreossi16}. Computational scalability varies among these different methods, with the most scalable methods requiring that the system dynamics do not involve control and disturbance variables. The work in \cite{Nilsson2016} accounts for both control and disturbances, but is only applicable to linear systems. Methods that can account for general nonlinear systems such as \cite{Althoff2014a} also sometimes represent sets using simple shapes such as polytopes, potentially sacrificing representation fidelity in favor of the other aspects mentioned earlier. Hamilton-Jacobi (HJ) formulations \cite{Barron90, Mitchell05, Margellos11, Bokanowski11} excel in handling general nonlinear dynamics, control and disturbance variables, and flexible set representations via a grid-based approach; however, these methods are the least computationally scalable. Still other methods make a variety of other assumptions to make desirable trade offs \cite{Darbon16, Hafner2009, Coogan2015}. In addition, under some special scenarios, it may be possible to obtain small computational benefits while minimizing trade offs in other axes of consideration by exploiting system structure \cite{Mitchell11, Fisac15, Mitchell03, Chen2016b, Kaynama2009, Kaynama2013}.}



In this paper, we present a system decomposition method for computing BRSs and BRTs of a class of nonlinear systems. Our method drastically reduces dimensionality without making any other trade offs. Our method first computes BRSs for lower-dimensional subsystems, and then reconstructs the full-dimensional BRS without incurring additional approximation errors other than those arising from the lower-dimensional computations. Crucially, the subsystems can be coupled through common states, controls, and disturbances. The treatment of this coupling distinguishes our method from others which consider completely decoupled subsystems, potentially obtained through transformations \cite{Callier1991, Sastry1999}. Since BRTs are also of great interest in many situations, we prove conditions under which BRTs can also be decomposed. 

\MCnote{The theory we present in this paper is compatible with any methods that compute BRSs and BRTs, such as \cite{Chen2013, Althoff2015, Frehse2011, Mitchell11, Chen2016b} and others mentioned earlier. In addition, when different decomposition methods are combined together, even more dimensionality reduction can be achieved.} This paper will be presented as follows:

\begin{itemize}
\item In Sections \ref{sec:background} and \ref{sec:formulation} we introduce the basic concept of reachability, and all the definitions needed for our proposed decomposition technique.
\item In Sections \ref{sec:sc} and \ref{sec:decoupled} we present our theoretical results related to decomposing BRSs for systems involving a control variable, but \textit{not} involving a disturbance variable.
\item In Section \ref{sec:set_to_tube} we show how BRTs can be decomposed.
\item In Section \ref{sec:hd} we demonstrate our decomposition method on high-dimensional systems.
\item In Section \ref{sec:dstb} we discuss how the presence of disturbances affects the above theoretical results.
\item We will also present numerical results obtained through the Hamilton-Jacobi (HJ) reachability formulation in \cite{Mitchell05} throughout the paper to validate our theory.
\end{itemize} 

\section{Background \label{sec:background}}
There are various formulations for computing the BRS and BRT when the system dimensionality is low. In this section, we give the basic mathematical problem setup to provide a foundation on which we build the new proposed theory.

\subsection{System Dynamics}
Let $\state\in \R^n$ be the system state, which evolves according to the ordinary differential equation (ODE)
\begin{equation}
\begin{aligned}
\label{eq:fdyn}
\frac{d\state(s)}{ds} = \dot\state(s) = \fdyn(\state(s), \ctrl(s)), s \in [t, 0], \ctrl(s) \in \cset
\end{aligned}
\end{equation}

In general, the theory we present is applicable when some states are periodic dimensions (such as angles), but for simplicity we will consider $\R^n$. The control is denoted by $\ctrl(s)$, with the control function $\ctrl(\cdot) \in \cfset$ being drawn from the set of measurable functions.


The control signal $\ctrl \in \cset \subset \mathbb{R}^{n_u}$ is compact and $t < 0$. The system dynamics, or flow field, $\fdyn: \R^n \times \cset \rightarrow \R^n$ is assumed to be uniformly continuous, bounded, and Lipschitz continuous in\footnote{For the remainder of the paper, we will omit the notation ``$(s)$'' from variables such as $\state$ and $\ctrl$ when referring to function values.} $\state$ for fixed $\ctrl$. Therefore, given $\ctrl(\cdot) \in \cfset$, there exists a unique trajectory solving \eqref{eq:fdyn} \cite{Coddington55}. We will denote solutions, or trajectories of \eqref{eq:fdyn} starting from state $\state$ at time $t$ under control $\ctrl(\cdot)$ as $\traj(s; \state, t, \ctrl(\cdot)): [t, 0] \rightarrow \R^n$. $\traj$ satisfies \eqref{eq:fdyn} with an initial condition almost everywhere:
\begin{equation}
\label{eq:fdyn_traj}
\begin{aligned}
\frac{d}{ds}\traj(s; \state, t, \ctrl(\cdot)) &= \fdyn(\traj(s; \state, t, \ctrl(\cdot)), \ctrl(s)) \\
\traj(t; \state, t, \ctrl(\cdot)) &= \state
\end{aligned}
\end{equation}

Since the dynamics \eqref{eq:fdyn} is time-invariant, the time variables in trajectories can also be shifted by some constant\footnote{In this case, it is implicit that the control function $\ctrl(\cdot)$ is also time-shifted by the same amount $\tau$.} $\tau$:
\begin{equation}
\label{eq:time-invariant}
\traj(s; \state, t, \ctrl(\cdot)) = \traj(s+\tau; \state, t+\tau, \ctrl(\cdot)), \forall \state \in \R^n
\end{equation}

\subsection{Backward Reachable Sets and Tubes}
\label{sec:RSRT}
We consider two different definitions of the BRS and two different definitions of the BRT. 

Intuitively, a BRS represents the set of states $\state\in\R^n$ from which the system can be driven into some set $\targetset \subseteq \R^n$ \textit{at the end} of a time horizon of duration $|t|$. We call $\targetset$ the ``target set''. First we define the ``Maximal BRS''; in this case the system seeks to enter $\targetset$ using some control function. We can think of $\targetset$ as a set of goal states. The Maximal BRS represents the set of states from which the system is guaranteed to reach $\targetset$. The second definition is for the ``Minimal BRS"; in this case the BRS is the set of states that will lead to $\targetset$ for all possible controls. Here we often consider $\targetset$ to be an unsafe set such as an obstacle. The Minimal BRS represents the set of states that leads to violation of safety requirements. Formally, the two definitions of BRSs are below\footnote{Sometimes in the literature, the argument of $\maxbrs$, $\minbrs$, $\maxbrt$, or $\minbrt$ is some non-negative number $\tau = -t$; however, for simplicity we will use the non-positive number $t$ to refer to the time horizon of the BRS and BRT.}:

\begin{defn}
\label{defn:rset_goal}
\textbf{Maximal BRS}.
\begin{equation*}
\maxbrs(t) = \{\state: \exists \ctrl(\cdot) \in \cfset, \traj(0; \state, t, \ctrl(\cdot)) \in \targetset \}
\end{equation*}
\end{defn}

\begin{defn}
\label{defn:rset_avoid}
\textbf{Minimal BRS}.
\begin{equation*}
\minbrs(t) = \{\state: \forall \ctrl(\cdot) \in \cfset, \traj(0; \state, t, \ctrl(\cdot)) \in \targetset \}
\end{equation*}
\end{defn}

\begin{figure}
	\centering
	\includegraphics[width=0.5\columnwidth]{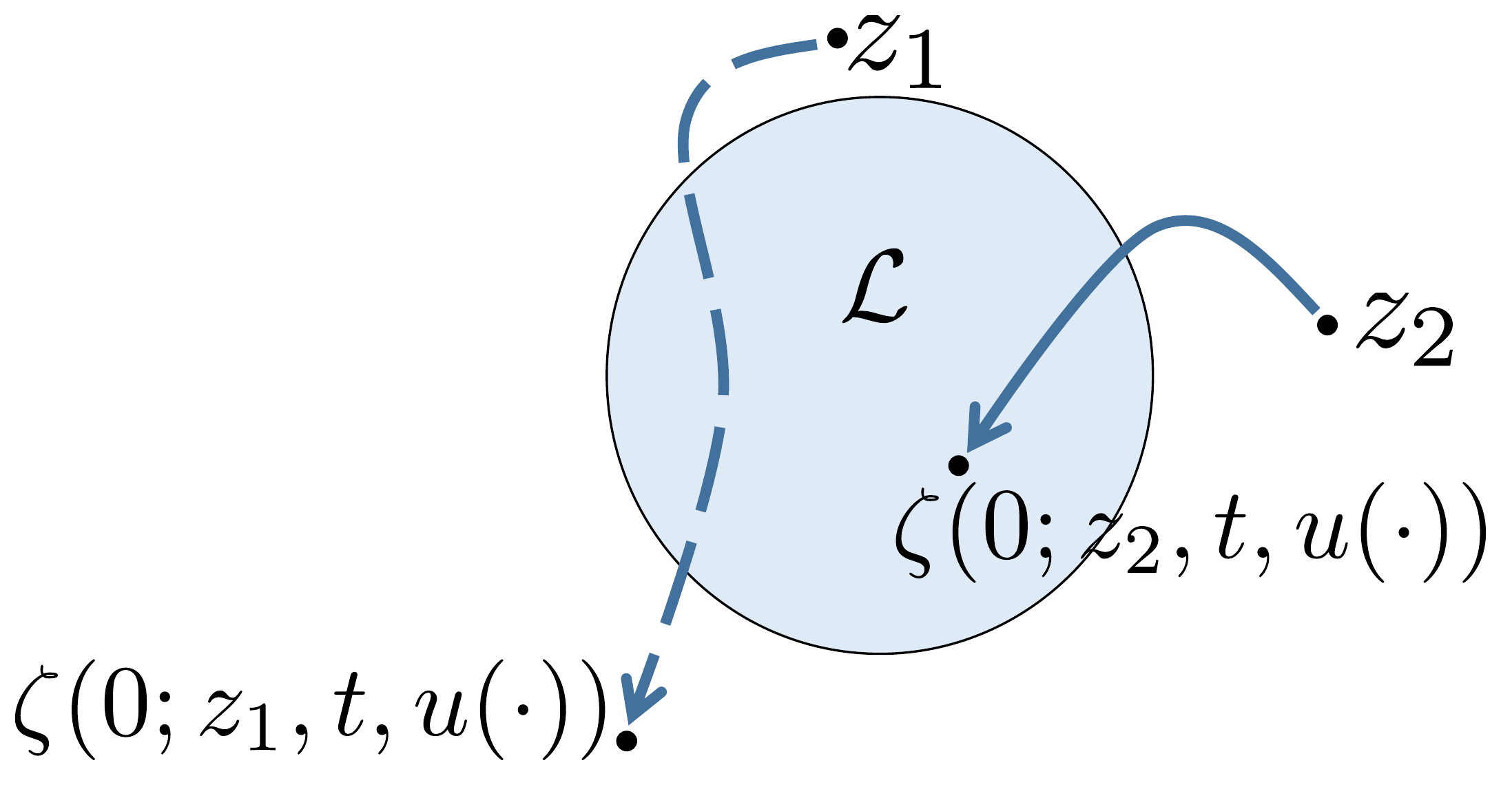}
	\caption{The difference between a BRS and a BRT. The dashed trajectory starts at $\state_1$ and passes through $\targetset$ during the period $[t, 0]$, but exits $\targetset$ by the end of the time period. Therefore the $\state_1$ is in the BRT, but not in the BRS. The solid trajectory starting from $\state_2$ is in $\targetset$ at the end of the time period. Therefore, $\state_2$ is in both the BRS and the BRT.}
	\label{fig:BRS_vs_BRT}
\end{figure}

While BRSs indicate whether a system can be driven into $\targetset$ at the end of a time horizon, BRTs indicate whether a system can be driven into $\targetset$ \textit{at some time} during the time horizon of duration $|t|$. Figure \ref{fig:BRS_vs_BRT} demonstrates the difference. BRTs are very important notions especially in safety-critical applications, in which we are interested in determining the ``Minimal BRT": the set of states that could lead to danger at some time within a specified time horizon. Formally, the two definitions of BRTs are as follows:

\begin{defn}
\label{defn:rtube_goal}
\textbf{Maximal BRT}.
\begin{equation*}
\maxbrt(t) = \{\state: \exists \ctrl(\cdot) \in \cfset, \exists s \in [t, 0], \traj(s; \state, t, \ctrl(\cdot)) \in \targetset \}
\end{equation*}
\end{defn}

\begin{defn}
\label{defn:rtube_avoid}
\textbf{Minimal BRT}.
\begin{equation*}
\minbrt(t) = \{\state: \forall \ctrl(\cdot) \in \cfset, \exists s \in [t, 0], \traj(s; \state, t, \ctrl(\cdot)) \in \targetset \}
\end{equation*}
\end{defn}

The terms ``maximal'' and ``minimal'' refer to the role of the optimal control \cite{Mitchell07b}. In the maximal (or minimal) case, the control causes the BRS or BRT to contain as many (or few) states as possible -- to have maximal (or minimal) size. 

\MCnote{While BRSs and BRTs indicate sets of states of interest, from a practical implementation perspective controller synthesis based on the reachable sets is extremely important. Much of the prior work on reachable set computations also include controller synthesis, which is usually done by casting the reachability problem as an optimal control or optimization problem, often with a functional representation of BRSs and BRTs. The controller is given as decision variables in the optimization \cite{Varaiya67, Barron90, Mitchell05,Bokanowski11, Majumdar13, Zhao2016}. We will not delve into the details of controller synthesis, since the theory we present in this paper is agnostic to these details.}

\section{Problem Formulation \label{sec:formulation}}
In this paper, we seek to obtain the BRSs and BRTs in Definitions \ref{defn:rset_goal} to \ref{defn:rtube_avoid} via computations in lower-dimensional subspaces under the assumption that the system \eqref{eq:fdyn} can be decomposed into self-contained subsystems (SCS) \eqref{eq:scdyn}. Such a decomposition is common, since many systems involve components that are loosely coupled. In particular, the evolution of position variables in vehicle dynamics is often weakly coupled though other variables such as heading. 

\subsection{Definitions}
\subsubsection{\textbf{Subsystem Dynamics}}
Let the state $\state \in \R^n$ be partitioned as $\state = (\spart_1, \spart_2, \spart_c)$, with $\spart_1 \in \R^{n_1}, \spart_2 \in \R^{n_2}, \spart_c \in \R^{n_c}, n_1, n_2 > 0, n_c \ge 0$, $n_1 + n_2 + n_c = n$. Note that $n_c$ could be zero. We call $\spart_1, \spart_2, \spart_c$ ``state partitions'' of the system. Intuitively, $\spart_1$ and $\spart_2,$ are states belonging to subsystems 1 and 2, respectively, and $\spart_c$ states belong to both subsystems.

Under the above notation, the system dynamics \eqref{eq:fdyn} become

\begin{equation}
\label{eq:expanded_dyn}
\begin{aligned}
\dot\spart_1 = \fdyn_1(\spart_1, \spart_2, \spart_c, \ctrl) \\
\dot\spart_2 = \fdyn_2(\spart_1, \spart_2, \spart_c, \ctrl) \\
\dot\spart_c = \fdyn_c(\spart_1, \spart_2, \spart_c, \ctrl) \\
\end{aligned}
\end{equation}

In general, depending on how the dynamics $\fdyn$ depend on $\ctrl$, some state partitions may be independent of the control. 

We group these states into subsystems by defining the SCS states $\scstate_1  = (\spart_1, \spart_c) \in \R^{n_1+n_c}$ and $\scstate_2 = (\spart_2, \spart_c) \in \R^{n_2+n_c}$, where $\scstate_1$ and $\scstate_2$ in general share the ``common'' states in $\spart_c$. Note that our theory is applicable to any finite number of subsystems defined in the analogous way, with $\scstate_i = (\spart_i, \spart_c)$; however, without loss of generality (WLOG), we assume that there are just two subsystems. 


\begin{defn}
\textbf{Self-contained subsystem}. Consider the following special case of \eqref{eq:expanded_dyn}:

\begin{equation}
\label{eq:scdyn}
\begin{aligned}
\dot\spart_1 &= \fdyn_1(\spart_1, \spart_c, \ctrl) \\
\dot\spart_2 &= \fdyn_2(\spart_2, \spart_c, \ctrl) \\
\dot\spart_c &= \fdyn_c(\spart_c, \ctrl) \\
\end{aligned}
\end{equation}

\end{defn}

We call each of the subsystems with states defined as $\scstate_i = (\spart_i, \spart_c)$ a ``self-contained subsystem'' (SCS), or just ``subsystem'' for short. Intuitively \eqref{eq:scdyn} means that the evolution of each subsystem depends only on the subsystem states: $\dot\scstate_i$ depends only on $\scstate_i = (\spart_i, \spart_c)$. Explicitly, the dynamics of the two subsystems are as follows:

\begin{equation*}
\begin{aligned}
&\dot\spart_1 = \fdyn_1(\spart_1, \spart_c, \ctrl) &\qquad \dot\spart_2 = \fdyn_2(\spart_2, \spart_c, \ctrl) \\
&\dot\spart_c = \fdyn_c(\spart_c, \ctrl) &\qquad \dot\spart_c = \fdyn_c(\spart_c, \ctrl)\\
&\text{(Subsystem 1)} & \text{(Subsystem 2)}
\end{aligned}
\end{equation*}

Note that the two subsystems are coupled through the common state partition $\spart_c$ and control $\ctrl$. When the subsystems are coupled through $\ctrl$, we say that the subsystems have ``shared control''. 

An example of a system that can be decomposed into SCSs is the Dubins Car with constant speed $\vel$:

\begin{equation}
\label{eq:dubins_car}
\begin{aligned}
\left[ \begin{array}{c}
\dot\pos_x\\
\dot\pos_y \\
\dot\theta
\end{array} \right]
=
\left[
\begin{array}{c}
v \cos\theta \\
v \sin\theta \\
\omega
\end{array}\right],
\qquad \omega \in \cset
\end{aligned}
\end{equation}

\noindent with state $\state = (\pos_x, \pos_y, \theta)$ representing the $x$ position, $y$ position, and heading, and control $\ctrl = \omega$ representing the turn rate. The state partitions are simply the system states: $\spart_1 = \pos_x, \spart_2 = \pos_y, \spart_c = \theta$. The subsystem states $\scstate_i$ and the subsystem controls $\scctrl_i$ are
\begin{equation}
\label{eq:dubins_car_decomp}
\begin{aligned}
\dot{\scstate_1} = 
\left[ \begin{array}{c}
\dot{\spart_1}\\
\dot{\spart_c}
\end{array} \right]
=
\left[ \begin{array}{c}
\dot{\pos_x}\\
\dot{\theta}
\end{array} \right]
&=
\left[\begin{array}{c}
v \cos\theta\\
\omega \\
\end{array}\right]\\
\dot{\scstate_2} = 
\left[ \begin{array}{c}
\dot{\spart_2}\\
\dot{\spart_c}
\end{array} \right]
=
\left[ \begin{array}{c}
\dot{\pos_y} \\
\dot{\theta}
\end{array} \right]
&=
\left[\begin{array}{c}
v \sin\theta \\
\omega
\end{array}\right]\\
\ctrl = \omega
\end{aligned}
\end{equation}

\noindent where the overlapping state is $\theta$, and the subsystem controls and their shared component is the control $\ctrl$ itself. The control partitions $\cpart_1, \cpart_2$ do not exist, since the state partitions $\spart_1, \spart_2$ do not depend on the control. For more examples of systems decomposed into SCSs, see \eqref{eq:quad6D}, \eqref{eq:quad6Dsc} and other numerical examples in this paper.

Although there may be common or overlapping states in $\scstate_1$ and $\scstate_2$, the evolution of each subsystem does not depend on the other explicitly. In fact, if we for example entirely ignore the subsystem $\scstate_2$, the evolution of the subsystem $\scstate_1$ is well-defined and can be considered a full system on its own; hence, each subsystem is self-contained.

\subsubsection{\textbf{Projection Operators}}
\begin{figure}[H]
	\centering
	\includegraphics[width=0.9\columnwidth, trim={0cm 8cm 1cm 8cm}, clip]{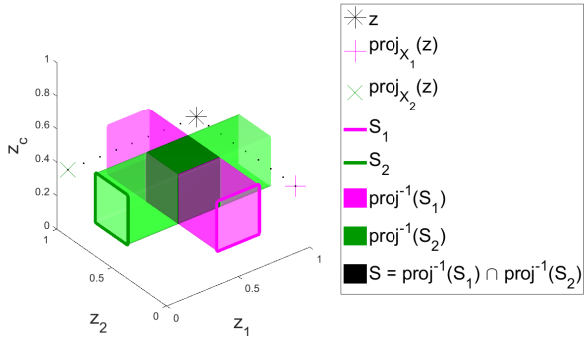}
	\caption{This figure shows the back-projection of sets in the $\state_1$-$\state_c$ plane $S_1$ and the $\state_2$-$\state_c$ plane $(S_2)$ to the 3D space to form the intersection shown as the black cube $(S)$. The figure also shows projection of a point $\state$ onto the lower-dimensional subspaces in the $\state_1$-$\state_c$ and $\state_2$-$\state_c$ planes.}
	\label{fig:proj}
\end{figure}

For the projection operators, it will be helpful to refer to Fig. \ref{fig:proj}. Define the projection of a state $\state = (\spart_1, \spart_2, \spart_c)$ onto a subsystem state space $\R^{n_i+n_c}$ as
\begin{equation}
\label{eq:project_pt}
\proj_i(\state) = \scstate_i = (\spart_i, \spart_c)
\end{equation}

This projects a point in the full dimensional state space onto a point in the subsystem state space. Also define the back-projection operator to be

\begin{equation}
\label{eq:backproject_pt}
\bp(\scstate_i) = \{\state \in \zset: (\spart_i, \spart_c) = \scstate_i\}
\end{equation}

This back-projection lifts a point from the subsystem state space to a set in the full dimensional state space. We will also need the ability to apply the back-projection operator on subsystems set to full dimensional sets. In this case, we overload the back-projection operator:

\begin{equation}
\label{eq:backproject_set}
\bp(\genset_i) = \{\state \in \zset: \exists \scstate_i \in \genset_i, (\spart_i, \spart_c) = \scstate_i\}
\end{equation}


\subsubsection{\textbf{Subsystem Trajectories}}
Since each subsystem in \eqref{eq:scdyn} is self-contained, we can denote the subsystem trajectories $\sctraj_i(s; \scstate_i, t, \ctrl(\cdot))$. When needed, we will write the subsystem trajectories more explicitly in terms of the state partitions as $\sctraj_i(s; \spart_i, \spart_c, t, \ctrl(\cdot))$. The subsystem trajectories satisfy the subsystem dynamics and initial condition:

\begin{equation}
\label{eq:scdyn_traj}
\begin{aligned}
\frac{d}{ds}\sctraj_i(s; \scstate_i, t, \ctrl(\cdot)) &= \scdyn_i(\sctraj_i(s; \scstate_i, t, \ctrl(s)) \\
\sctraj_i(t; \scstate_i, t, \ctrl(\cdot)) &= \scstate_i
\end{aligned}
\end{equation}

\noindent where $\scdyn_i(\scstate_i, \ctrl) = (\fdyn_i(\spart_i, \spart_c, \ctrl), \fdyn_c(\spart_c, \ctrl))$, and the full system trajectory and subsystem trajectories are simply related to each other via the projection operator:
\begin{equation}
\label{eq:proj_traj}
\proj_i\big(\traj(s; \state, t, \ctrl(\cdot))\big) = \sctraj_i(s; \scstate_i, t, \ctrl(\cdot))
\end{equation}

\noindent where $\scstate_i = \proj_i(\state)$.

\subsection{Goals of This Paper}
We assume that the full system target set $\targetset$ can be written in terms of the subsystem target sets $\sctarget{1} \subseteq \xset_1, \sctarget{2} \subseteq \xset_2$ in one of the following ways:

\begin{equation}
\label{eq:target_intersect}
\targetset = \bp(\sctarget{1}) \cap \bp(\sctarget{2})
\end{equation}

\noindent where the full target set is the intersection of the back-projections of subsystem target sets, or

\begin{equation}
\label{eq:target_union}
\targetset = \bp(\sctarget{1}) \cup \bp(\sctarget{2})
\end{equation}

\noindent where the full target set is the union of the back-projections of subsystem target sets. Fig. \ref{fig:proj} helps provide intuition for these concepts: applying (\ref{eq:target_intersect}) to $\genset_1$ and $\genset_2$ results in the black cube.  Applying (\ref{eq:target_union}) would result in the cross-shaped set encompassing both $\bp(\genset_1)$ and $\bp(\genset_2)$.

In practice, this is not a strong assumption, since $\sctarget{1}$ and $\sctarget{2}$ share the common variables $\spart_c$. Relatively complex shapes, for example those in Fig. \ref{fig:dubins_angles} and \ref{fig:dubins_un_set}, can be represented by an intersection or union of back-projections of lower-dimensional sets that share common variables. In addition, such an assumption is reasonable since the full system target set should at least be representable in some way in the lower-dimensional spaces. 

Next, we define the subsystem BRSs $\maxbrs_i, \minbrs_i$ the same way as in Definitions \ref{defn:rset_goal} and \ref{defn:rset_avoid}, but with the subsystems in \eqref{eq:scdyn} and subsystem target sets $\sctarget{i}, i=1,2$, respectively:

\begin{equation}
\begin{aligned}
\label{eq:sc_rset}
\maxbrs_i(t) &= \{\scstate_i: \exists \ctrl(\cdot), \sctraj_i(0; \scstate_i, t, \ctrl(\cdot)) \in \sctarget{i} \} \\
\minbrs_i(t) &= \{\scstate_i: \forall \ctrl(\cdot), \sctraj_i(0; \scstate_i, t, \ctrl(\cdot)) \in \sctarget{i} \}
\end{aligned}
\end{equation}

Subsystem BRTs are defined analogously:
\begin{equation}
\begin{aligned}
\label{eq:sc_rtube}
\maxbrt_i(t) = \{\scstate_i: \exists \ctrl(\cdot), \exists s \in [t, 0], \sctraj_i(s; \scstate_i, t, \ctrl(\cdot)) \in \sctarget{i} \} \\
\minbrt_i(t) = \{\scstate_i: \forall \ctrl(\cdot), \exists s \in [t, 0], \sctraj_i(s; \scstate_i, t, \ctrl(\cdot)) \in \sctarget{i} \}
\end{aligned}
\end{equation}


Given a system in the form of \eqref{eq:scdyn} with target set that can be represented by \eqref{eq:target_intersect} or \eqref{eq:target_union}, our goals are as follows.



\begin{itemize}
\item \textbf{Decomposition of BRSs.} First, we would like to compute full-dimensional BRSs by performing computations in lower-dimensional subspaces. Specifically, we would like to first compute the subsystem BRSs $\maxbrs_i(t)$ or $\minbrs_i(t)$, and then reconstruct the full system BRS $\maxbrs(t)$ or $\minbrs(t)$. This process greatly reduces computation burden by decomposing the full system into two lower-dimensional subsystems. Formally, we would like to investigate the situations in which the following four cases is true:

\begin{equation}
\begin{aligned}
\label{eq:recon}
\text{\eqref{eq:target_intersect}} &\Rightarrow \maxbrs(t) = \bp(\maxbrs_1(t)) \cap \bp(\maxbrs_2(t)) \\
\text{\eqref{eq:target_intersect}} &\Rightarrow \minbrs(t) = \bp(\minbrs_1(t)) \cap \bp(\minbrs_2(t)) \\
\text{\eqref{eq:target_union}} &\Rightarrow \maxbrs(t) = \bp(\maxbrs_1(t)) \cup \bp(\maxbrs_2(t)) \\
\text{\eqref{eq:target_union}} &\Rightarrow \minbrs(t) = \bp(\minbrs_1(t)) \cup \bp(\minbrs_2(t))
\end{aligned}
\end{equation}

Results related to BRSs are outlined for SCSs in Theorems \ref{thm:sc_reach_un} and \ref{thm:sc_reach_in}. In the case that the subsystem controls do not share any components, Propositions \ref{prop:decoupled_reach_in} and \ref{prop:decoupled_reach_un} state stronger results.

\item \textbf{Decomposition of BRTs}. BRTs are useful since they provide guarantees over a time horizon as opposed to at a particular time. However, often BRTs cannot be decomposed the same way as BRSs. Therefore, our second goal is to propose how BRTs can be decomposed. These results are stated in Propositions \ref{prop:BRT_Union} and \ref{prop:sets2tube_goal}, and Theorem \ref{thm:sets2tube_avoid}.
\item \textbf{Treatment of disturbances}. Finally, we investigate how the above theoretical results change in the presence of disturbances. In Section \ref{sec:dstb}, we will show that slightly conservative BRSs and BRTs can still be obtained using our decomposition technique.
\end{itemize}

Tables \ref{table:summary_BRS} and \ref{table:summary_BRT} summarize our theoretical results and where details of each result can be found.

\begin{table*}[]
	\centering
	\caption{Backward Reachable Set Decomposition}
	\label{table:summary_BRS}
	\SHnote{
	\begin{tabular}{|l|l|l|l|l|l|l|l|l|}
		\hline
		\textbf{Section}            & \multicolumn{2}{c|}{\textbf{\ref{sec:sc}}}   & \multicolumn{2}{c|}{\textbf{\ref{sec:decoupled}}}    & \multicolumn{2}{c|}{\textbf{\ref{sec:sc_dstb}}} & \multicolumn{2}{c|}{\textbf{\ref{sec:decoupled_dstb}}} \\ \hline
		\textbf{Shared Controls}    & \multicolumn{2}{c|}{\textbf{Yes}}  & \multicolumn{2}{c|}{\textbf{No}}   & \multicolumn{2}{c|}{\textbf{Yes}}    & \multicolumn{2}{c|}{\textbf{No}}     \\ \hline
		\textbf{Shared Disturbance} & \multicolumn{2}{c|}{\textbf{No}}   & \multicolumn{2}{c|}{\textbf{No}}   & \multicolumn{2}{c|}{\textbf{Yes}}    & \multicolumn{2}{c|}{\textbf{Yes}}    \\ \hline
		\textbf{Target}             & \textbf{Intersection} & \textbf{Union}      & \textbf{Intersection} & \textbf{Union}      & \textbf{Intersection}  & \textbf{Union}       & \textbf{Intersection}  & \textbf{Union}       \\ \hline
		Recover Max. BRS?  & No           & Yes, exact & Yes, exact   & Yes, exact & No            & Yes, consrv & Yes, consrv   & Yes, consrv          \\ \hline
		Recover Min. BRS?  & Yes, exact   & No         & Yes, exact   & Yes, exact & Yes, consrv   & No          & Yes, consrv            & Yes, consrv \\ \hline
		Locations \& Equation(s)        & Thm \ref{thm:sc_reach_in}, (\ref{eq:targetset_in})  & Thm \ref{thm:sc_reach_un}, (\ref{eq:targetset_un})         & \begin{tabular}[c]{@{}l@{}}Prop \ref{prop:decoupled_reach_in}, (\ref{eq:targetset_in_decoupled})\\ Thm \ref{thm:sc_reach_in}, (\ref{eq:targetset_in})\end{tabular} & \begin{tabular}[c]{@{}l@{}}Thm \ref{thm:sc_reach_un}, (\ref{eq:targetset_un})\\ Prop \ref{prop:decoupled_reach_un}, (\ref{eq:targetset_un_decoupled})\end{tabular}          & Cor \ref{cor:minBRS_Intersection_SCS_Disturbance}, (\ref{eq:minBRS_intersection_disturbance}) & Cor \ref{cor:maxBRS_Union_SCS_Disturbance},  (\ref{eq:maxBRS_union_disturbance}) & Cor \ref{cor:targetset_in_decoupled_with_dstb}, (\ref{eq:targetset_in_decoupled_with_dstb}) & Cor \ref{cor:targetset_un_decoupled_with_dstb}, (\ref{eq:targetset_un_decoupled_with_dstb})           \\ \hline
	\end{tabular}
	\bigskip
	\begin{flushleft}
	 Summary of possible decompositions of the BRS, whether they are possible, and if so whether they are exact or conservative. Exact means that no additional approximation errors are introduced. Note that in the cases marked ``no" for shared control (or shared disturbance), the results hold for both decoupled control (or disturbance) and for no control (or disturbance). All cases shown are for scenarios with shared states, with the shared states being $\spart_c$ in \eqref{eq:scdyn}; in the case that there are no shared states this becomes a straightforward decoupled system.
	\end{flushleft}
}
\end{table*}

\begin{table*}[]
	\centering
	\caption{BRT Results for Reconstruction from Tubes}
\label{table:BRT_Tubes_Summary}
	\SHnote{\begin{tabular}{|l|l|l|l|l|c|l|l|l|}
		\hline
		\textbf{Section}            & \multicolumn{4}{c|}{\textbf{\ref{subsec:tube2tube}}}                               & \multicolumn{4}{c|}{\textbf{\ref{sec:set2tube_dstb}}}                                                    \\ \hline
		\textbf{Shared Controls}    & \multicolumn{2}{c|}{\textbf{Yes}}  & \multicolumn{2}{c|}{\textbf{No}}   & \multicolumn{2}{c|}{Yes}                         & \multicolumn{2}{c|}{\textbf{No}}     \\ \hline
		\textbf{Shared Disturbance} & \multicolumn{2}{c|}{\textbf{No}}   & \multicolumn{2}{c|}{\textbf{No}}   & \multicolumn{2}{c|}{Yes}                         & \multicolumn{2}{c|}{\textbf{Yes}}    \\ \hline
		\textbf{Target}             & \textbf{Intersection} & \textbf{Union}      & \textbf{Intersection} & \textbf{Union}      & \multicolumn{1}{l|}{\textbf{Intersection}} & \textbf{Union}        & \textbf{Intersection} & Union        \\ \hline
		Recover Max. BRT?  & No           & Yes, exact & No           & Yes, exact & No                                & Yes, conserv & No           & Yes, conserv \\ \hline
		Recover Min. BRT?  & No           & No         & No           & Yes, exact & No                                & No & No           & Yes, conserv \\ \hline
		Equation(s)        & N/A          & Prop \ref{prop:BRT_Union}, (\ref{eq:maxBRT_union})          & N/A          & \begin{tabular}[c]{@{}l@{}}Prop \ref{prop:BRT_Union}, (\ref{eq:maxBRT_union})\\ Prop \ref{prop:BRT_Union}, (\ref{eq:minBRT_union})\end{tabular}          & N/A                               & Cor \ref{cor:maxBRT_union_disturbance}, (\ref{eq:maxBRT_union_disturbance})            & N/A          & \begin{tabular}[c]{@{}l@{}}Cor \ref{cor:maxBRT_union_disturbance}, (\ref{eq:maxBRT_union_disturbance})\\ Cor \ref{cor:minBRT_union_disturbance}, (\ref{eq:minBRT_union_disturbance})\end{tabular}           \\ \hline
	\end{tabular}}
\end{table*}

\begin{table*}[]
	\centering
	\caption{BRT Results for Reconstruction from Sets}
\label{table:BRT_Sets_Summary}
	\SHnote{\begin{tabular}{|l|l|l|c|l|c|l|c|l|}
		\hline
		\textbf{Section}            & \multicolumn{4}{c|}{\textbf{\ref{subsec:union_BRT}}}                                              & \multicolumn{4}{c|}{\textbf{\ref{sec:set2tube_dstb}}}                                                           \\ \hline
		\textbf{Shared Controls}    & \multicolumn{2}{c|}{\textbf{Yes}} & \multicolumn{2}{c|}{\textbf{No}}                   & \multicolumn{2}{c|}{\textbf{Yes}}                  & \multicolumn{2}{c|}{\textbf{No}}                   \\ \hline
		\textbf{Shared Disturbance} & \multicolumn{2}{c|}{\textbf{No}}  & \multicolumn{2}{c|}{\textbf{No}}                   & \multicolumn{2}{c|}{\textbf{Yes}}                  & \multicolumn{2}{c|}{\textbf{Yes}}                  \\ \hline
		\textbf{Target}             & \textbf{Intersection}   & \textbf{Union}   & \multicolumn{1}{l|}{\textbf{Intersection}} & \textbf{Union} & \multicolumn{1}{l|}{\textbf{Intersection}} & \textbf{Union} & \multicolumn{1}{l|}{\textbf{Intersection}} & \textbf{Union} \\ \hline
		Recover Max. BRT?  & \multicolumn{4}{c|}{Yes, exact}                                      & \multicolumn{4}{c|}{Yes, conserv}                                                     \\ \hline
		Recover Min. BRT?  & \multicolumn{4}{c|}{Yes, exact*}                                     & \multicolumn{4}{c|}{Yes, exact*}                                                      \\ \hline
		Equation(s)        & \multicolumn{4}{c|}{\begin{tabular}[c]{@{}l@{}}Prop \ref{prop:sets2tube_goal}, (\ref{eq:maxBRT_intersection})\\ Thm \ref{thm:sets2tube_avoid}, (\ref{eq:minBRT_intersection})\end{tabular}}                                             & \multicolumn{4}{c|}{\begin{tabular}[c]{@{}l@{}}Cor \ref{cor:maxBRT_intersection_disturbance}, (\ref{eq:maxBRT_intersection_disturbance})\\ Thm \ref{thm:sets2tube_avoid}, (\ref{eq:minBRT_intersection})\end{tabular}}                                                               \\ \hline
	\end{tabular}
	\bigskip\begin{flushleft}
  Summary of possible decompositions of the BRT, whether they are possible, and if so whether they are exact or conservative. Exact means that no additional approximation errors are introduced. Note that in the cases marked ``no" for shared control (or shared disturbance), the results hold for both decoupled control (or disturbance) and for no control (or disturbance). All cases shown are for scenarios with shared states, with the shared states being $\spart_c$ in \eqref{eq:scdyn}; in the case that there are no shared states this becomes a straightforward decoupled system. We assume that exact sets are available to compute those BRTs that require the union of sets.\\
  * the solution here can be found only if the minimum BRSs are non-empty for the entire time period.
\end{flushleft}}
\end{table*}

\section{Self-Contained Subsystems \label{sec:sc}}
Suppose the full system \eqref{eq:fdyn} can be decomposed into SCSs given in \eqref{eq:scdyn}. Then, the full-dimensional BRS can be reconstructed without incurring additional approximation errors from lower-dimensional BRSs in the situations stated in Theorem \ref{thm:sc_reach_un} and \ref{thm:sc_reach_in}.

\begin{rem}
  If $\targetset$ represents states the system aims to reach, then $\maxbrs(t)$ represents the set of states from which $\targetset$ can be reached. If the system goal states are the union of subsystem goal states, then it suffices for any subsystem to reach its subsystem goal states, regardless of any coupling that exists between the subsystems. Theorem \ref{thm:sc_reach_un} states this intuitive result.
\end{rem}

\begin{thm}
\label{thm:sc_reach_un}
Suppose that the full system in \eqref{eq:fdyn} can be decomposed into the form of \eqref{eq:scdyn}, then

\begin{equation}
\label{eq:targetset_un}
\begin{aligned}
&\targetset = \bp(\sctarget{1}) \cup \bp(\sctarget{2}) \\
&\Rightarrow \maxbrs(t) = \bp(\maxbrs_1(t)) \cup \bp(\maxbrs_2(t))
\end{aligned}
\end{equation}
\end{thm}

\begin{rem}
  If $\targetset$ represents the set of unsafe states, then $\minbrs(t)$ is the set of states from which the system will be driven into danger. Thus outside of $\minbrs(t)$, there exists a control for the system to avoid the unsafe states. For the system to avoid $\targetset$, it suffices to avoid the unsafe states in either subsystem, regardless of any coupling that exists between the subsystems. Theorem \ref{thm:sc_reach_in} formally states this intuitive result.
\end{rem}

\begin{thm}
  \label{thm:sc_reach_in}
  Suppose that the full system in \eqref{eq:fdyn} can be decomposed into the form of \eqref{eq:scdyn}, then
  
  \begin{equation}
  \label{eq:targetset_in}
  \begin{aligned}
  &\targetset = \bp(\sctarget{1}) \cap \bp(\sctarget{2}) \\
  &\Rightarrow \minbrs(t) = \bp(\minbrs_1(t)) \cap \bp(\minbrs_2(t))
  \end{aligned}
  \end{equation}
\end{thm}

To prove the theorems, we need some intermediate results.

\begin{lem}
\label{lem:proj_basic}
Let $\bar\state\in\zset, \bar\scstate_i = \proj_i(\bar\state), \genset_i \subseteq \xset_i$. Then,

\begin{equation}
\bar\scstate_i \in \genset_i \Leftrightarrow \bar\state \in \bp(\genset_i)
\end{equation}
\end{lem}

\textit{Proof of Lemma \ref{lem:proj_basic}}: Forward direction: Suppose $\bar\scstate_i \in \genset_i$, then trivially $\exists \scstate_i \in \genset_i, \proj_i(\bar\state) = \scstate_i$. By the definition of back-projection in \eqref{eq:backproject_set}, we have $\bar\state \in \bp(\genset_i)$.

Backward direction: Suppose $\bar\state \in \bp(\genset_i)$, then by \eqref{eq:backproject_set} we have $\exists \scstate_i \in \genset_i, \proj_i(\bar\state) = \scstate_i$. Denote such an $\scstate_i$ to be $\hat\scstate_i$, and suppose $\bar\scstate_i \notin \genset_i$. Then, we have $\hat\scstate_i \neq \bar\scstate_i$, a contradiction, since $\bar\scstate_i = \proj_i(\bar\state) = \hat\scstate_i$. \hfill\IEEEQEDhere

\begin{cor}
  \label{cor:proj_un}
  If $\genset = \bp(\genset_1) \cup \bp(\genset_2)$, then
  \begin{equation*}
  \bar\state \in \genset \Leftrightarrow \bar\scstate_1\in\genset_1 \vee \bar\scstate_2\in\genset_2, \text{ where } \bar\scstate_i = \proj_i(\bar\state)
  \end{equation*}
\end{cor}

\begin{cor}
\label{cor:proj_in}
If $\genset = \bp(\genset_1) \cap \bp(\genset_2)$, then 
\begin{equation*}
\bar\state \in \genset \Leftrightarrow \bar\scstate_1 \in \genset_1 \wedge \bar\scstate_2 \in \genset_2, \text{ where } \bar\scstate_i = \proj_i(\bar\state)
\end{equation*}
\end{cor}

\subsection{Proof of Theorem \ref{thm:sc_reach_un}}
We will prove the following equivalent statement:
\begin{equation}
\label{eq:sc_proof_un}
\bar\state \in \maxbrs(t) \Leftrightarrow \bar\state \in \bp(\maxbrs_1(t)) \cup \bp(\maxbrs_2(t))
\end{equation}

Consider the relationship between the full system trajectory and subsystem trajectory in \eqref{eq:proj_traj}. Define $\bar\scstate_i = \proj_i(\bar\state)$ and $\sctraj_i(0; \bar\scstate_i, t, \ctrl(\cdot)) = \proj_i(\traj(0; \bar\state, t, \ctrl(\cdot)))$.

We first prove the backward direction. By Corollary \ref{cor:proj_un}, \eqref{eq:sc_proof_un} is equivalent to

\begin{equation}
\label{eq:subsystem_in_reachset}
\bar\scstate_1 \in \maxbrs_1(t)\vee\bar\scstate_2 \in \maxbrs_2(t)
\end{equation}

WLOG, assume $\bar\scstate_1 \in \maxbrs_1(t)$. By the subsystem BRS definition in \eqref{eq:sc_rset}, this is equivalent to

\begin{equation}
\label{eq:pick_system_1}
\begin{array}{rc}
\exists \ctrl(\cdot), \sctraj_1(0; \bar\scstate_1, t, \ctrl(\cdot)) \in \sctarget{1}\\
\end{array}
\end{equation}

By Lemma \ref{lem:proj_basic}, we equivalently have $\bar\state \in \bp(\maxbrs_1(t))$. This proves the backward direction. 

For the forward direction, we begin with $\bar\state \in \maxbrs(t)$, which by Definition \ref{defn:rset_goal} is equivalent to $
\exists \ctrl(\cdot), \traj(0; \bar\state, t, \ctrl(\cdot)) \in \targetset$. By Corollary \ref{cor:proj_un}, we then have

\begin{equation}
\label{eq:subsystem_reach_target}
\exists \ctrl(\cdot), \sctraj_1(0; \bar\scstate_1, t, \ctrl(\cdot)) \in \sctarget{1} \vee \sctraj_2(0; \bar\scstate_2, t, \ctrl(\cdot)) \in \sctarget{2}
\end{equation}

Finally, distributing ``$\exists \ctrl(\cdot)$'' gives \eqref{eq:subsystem_in_reachset}. \hfill\IEEEQEDhere

\subsection{Proof of Theorem \ref{thm:sc_reach_in}}
We will prove the following equivalent statement:
\begin{equation}
\label{eq:proof_left}
\bar\state \notin \minbrs(t) \Leftrightarrow \bar\state \notin \bp(\minbrs_1(t)) \cap \bp(\minbrs_2(t)) \\
\end{equation}

The above statement is equivalent to
\begin{equation}
\bar\state \in \minbrs^c(t) \Leftrightarrow \bar\state \in \big[\bp(\minbrs_1(t))\big]^c \cup \big[\bp(\minbrs_2(t))\big]^c
\end{equation}

By the Definition \ref{defn:rset_avoid} (minimal BRS), we have that $\bar\state \in \minbrs^c(t)$ is equivalent to $\exists \ctrl(\cdot) \in \cfset, \traj(0; \bar\state, t, \ctrl(\cdot)) \in \targetset^c$. Also, 

\begin{equation*}
\bar\state \in \big[\bp(\minbrs_1(t))\big]^c \cup \big[\bp(\minbrs_2(t))\big]^c
\end{equation*}

\noindent is equivalent to $\bar\scstate_1 \in \minbrs^c_1(t) \vee \bar\scstate_2 \in \minbrs^c_2(t)$.

From here, we can proceed in the same fashion as the proof of Theorem \ref{thm:sc_reach_un}, with ``$\maxbrs(t)$'' replaced with ``$\minbrs^c(t)$'', ``$\bp(\maxbrs_i(t))$'' replaced with ``$\big[\bp(\minbrs_i(t))\big]^c$'', and ``$\targetset$'', ``$\sctarget{i}$'' replaced with ``$\targetset^c$'', ``$\sctarget{i}^c$'', respectively. \hfill\IEEEQEDhere\\
The conditions for reconstruction the maximal BRS for an intersection of targets, as well as the minimal BRS for a union of targets, are more complicated and beyond the scope of this paper.  Table \ref{table:SCS_BRS} summarizes the results from this section.
\begin{table}[h]
	\centering
		\caption{BRS Results from Section \ref{sec:sc}}
		\label{table:SCS_BRS}
		\SHnote{\begin{tabular}{|l|l|l|}
			\hline
			\textbf{Shared Controls}       & \multicolumn{2}{c|}{\textbf{Yes}}   \\ \hline
			\textbf{Shared Disturbance}    & \multicolumn{2}{c|}{\textbf{No}}    \\ \hline
			\textbf{Target}                & \textbf{Intersection} & \textbf{Union}       \\ \hline
			Recover Max. BRS?     & No           & Yes, exact  \\ \hline
			Recover Min. BRS?     & Yes, exact   & No          \\ \hline
			Equation(s) & Thm \ref{thm:sc_reach_in}, (\ref{eq:targetset_in})  & Thm \ref{thm:sc_reach_un}, (\ref{eq:targetset_un}) \\ \hline
		\end{tabular}}
\end{table}

\subsection{Numerical Example: The Dubins Car}
The Dubins Car is a well-known system whose dynamics are given by \eqref{eq:dubins_car}. This system is only 3D, and its BRS can be tractably computed in the full-dimensional space, so we use it to compare the full formulation with our decomposition method. The Dubins Car dynamics can be decomposed according to \eqref{eq:dubins_car_decomp}. For this example, we computed the BRS from the target set representing positions near the origin in both the $\pos_x$ and $\pos_y$ dimensions:

\begin{equation}
\label{eq:dubins_target}
\targetset = \{(\pos_x, \pos_y, \theta): |\pos_x|, |\pos_y| \le 0.5\}
\end{equation}

\begin{figure}
	\centering
	\includegraphics[width=\columnwidth]{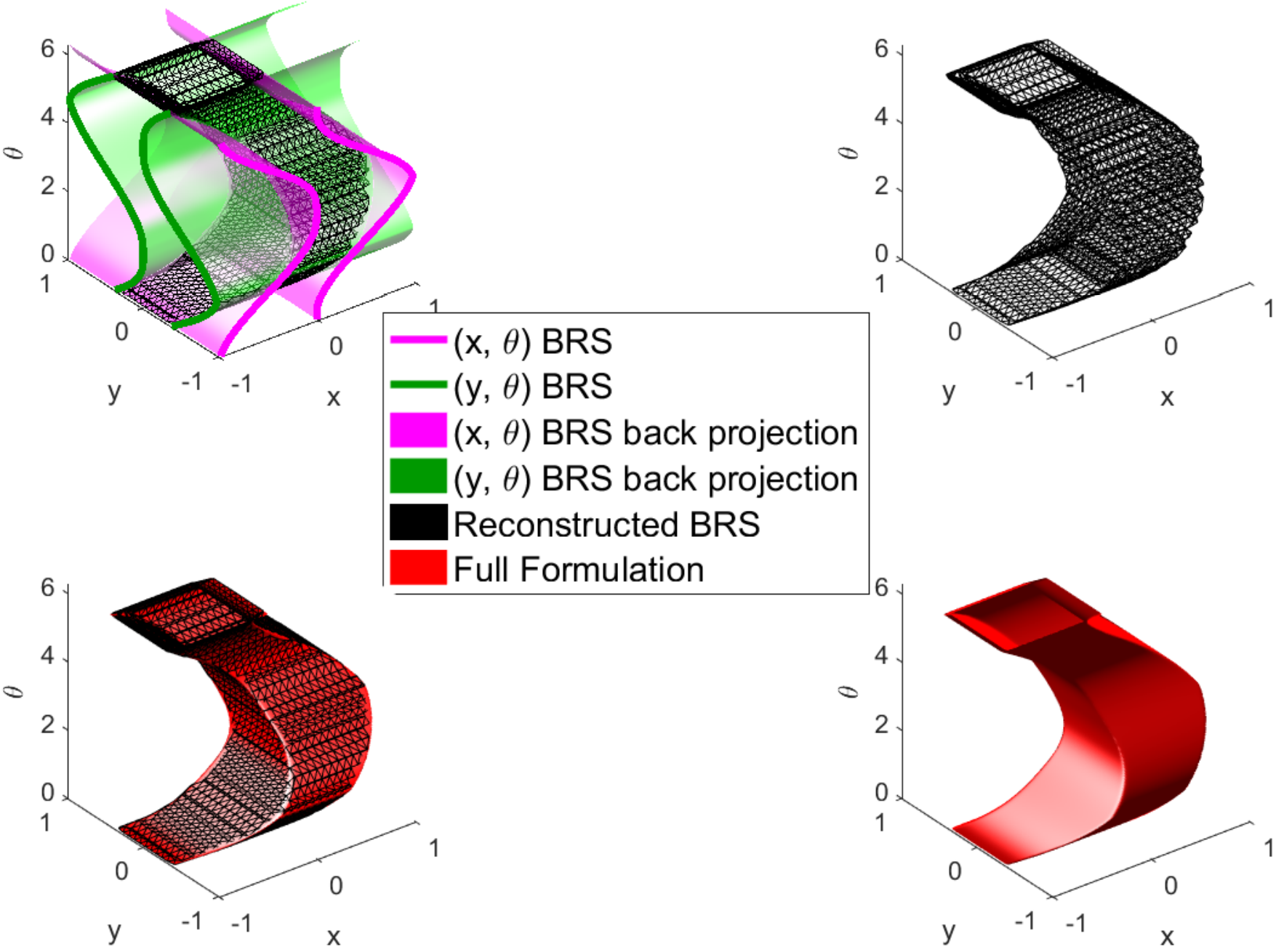}
	\caption{Comparison of the Dubins Car BRS $\minbrs(t=-0.5)$ computed using the full formulation and via decomposition. Left top: BRSs in the lower-dimensional subspaces and how they are combined to form the full-dimensional BRS. Top right: BRS computed via decomposition. Bottom left: BRSs computed using both methods, superimposed, showing that they are indistinguishable. Bottom right: BRS computed using the full formulation.}
	\label{fig:dubins_compare}
\end{figure}

Such a target set $\targetset$ can be used to model an obstacle that the vehicle must avoid. Given $\targetset$, the interpretation of the BRS $\minbrs(t)$ is the set of states from which a collision with the obstacle may occur after a duration of $|t|$. From $\targetset$, we computed the BRS $\minbrs(t)$ at $t=-0.5$. The resulting full formulation BRS is shown in Fig. \ref{fig:dubins_compare} as the red surface which appears in the bottom subplots. To compute the BRS using our decomposition method, we write the unsafe set $\targetset$ as

\begin{equation}
\label{eq:dubins_target_decomp}
    \begin{aligned}
    \sctarget{1} &= \{(\pos_x, \theta): |\pos_x| \le 0.5\}, \sctarget{2} = \{(\pos_y, \theta ): |\pos_y| \le 0.5\} \\
    \targetset &= \bp(\sctarget{1}) \cap \bp(\sctarget{2})
    \end{aligned}
\end{equation}

\begin{figure}
	\centering
	\includegraphics[width=0.9\columnwidth]{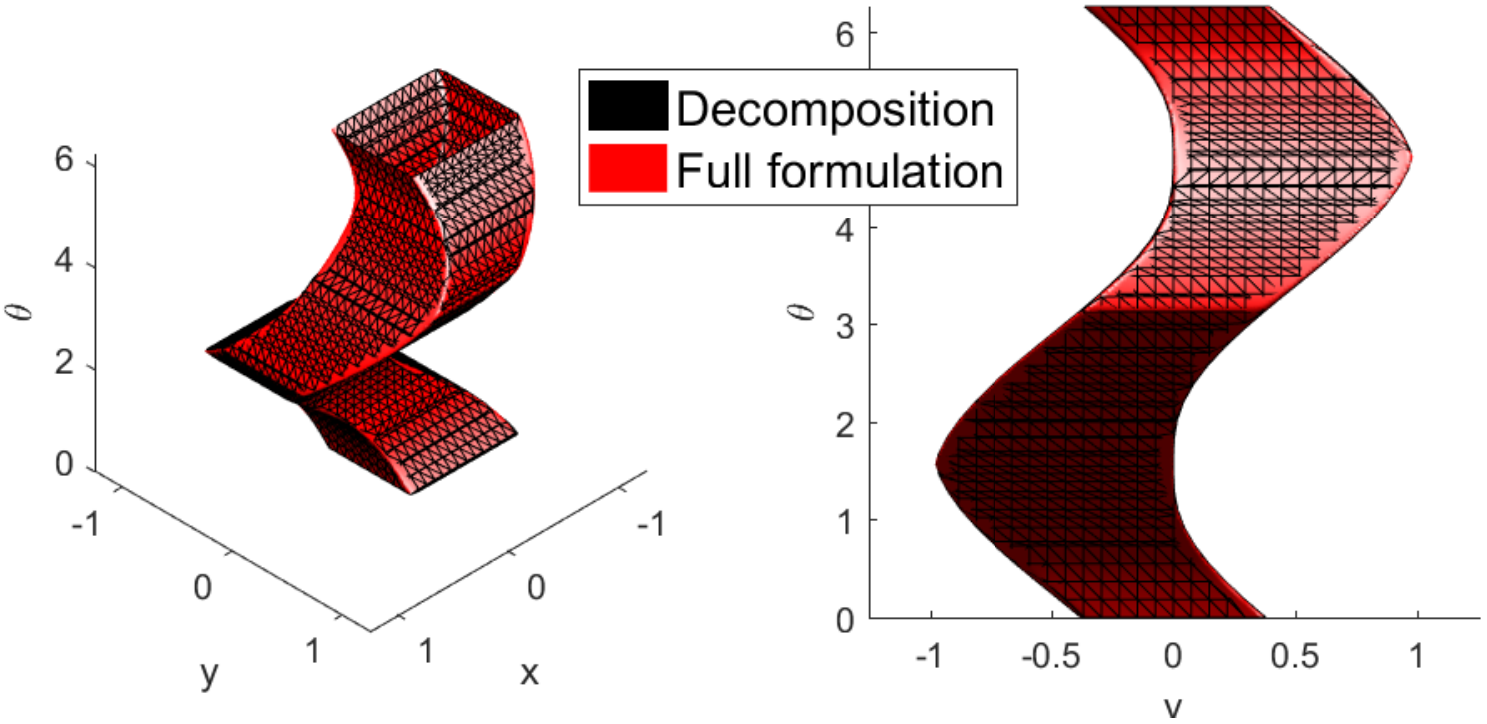}
	\caption{The Dubins Car BRS $\minbrs(t=-0.5)$ computed using the full formulation and via decomposition, other view angles.}
	\label{fig:dubins_angles}
\end{figure}

From $\sctarget{1}$ and $\sctarget{2}$, we computed the lower-dimensional BRSs $\minbrs_1(t)$ and $\minbrs_2(t)$, and then reconstructed the full-dimensional BRS $\minbrs(t)$ using Theorem \ref{thm:sc_reach_in}: $\minbrs(t) = \bp(\minbrs_1(t)) \cap \bp(\minbrs_2(t))$. The subsystem BRSs and their back-projections are shown in magenta and green in the top left subplot of Fig. \ref{fig:dubins_compare}. The reconstructed BRS is shown in the top left, top right, and bottom left subplots of Fig. \ref{fig:dubins_compare} (black mesh).

In the bottom left subplot of Fig. \ref{fig:dubins_compare}, we superimpose the full-dimensional BRS computed using the two methods. We show the comparison of the computation results viewed from two different angles in Fig. \ref{fig:dubins_angles}. The results are indistinguishable. 

Theorem \ref{thm:sc_reach_in} allows the computation to be performed in lower-dimensional subspaces, which is significantly faster. Another benefit of the decompositino method is that in the numerical methods for solving the HJ PDE, the amount of numerical dissipation increases with the number of state dimensions. Thus, computations in lower-dimensional subspaces lead to a slightly more accurate numerical solution.

\begin{figure}
	\centering
	\includegraphics[width=0.8\columnwidth, trim={0cm 0cm 0cm 5.25cm}, clip]{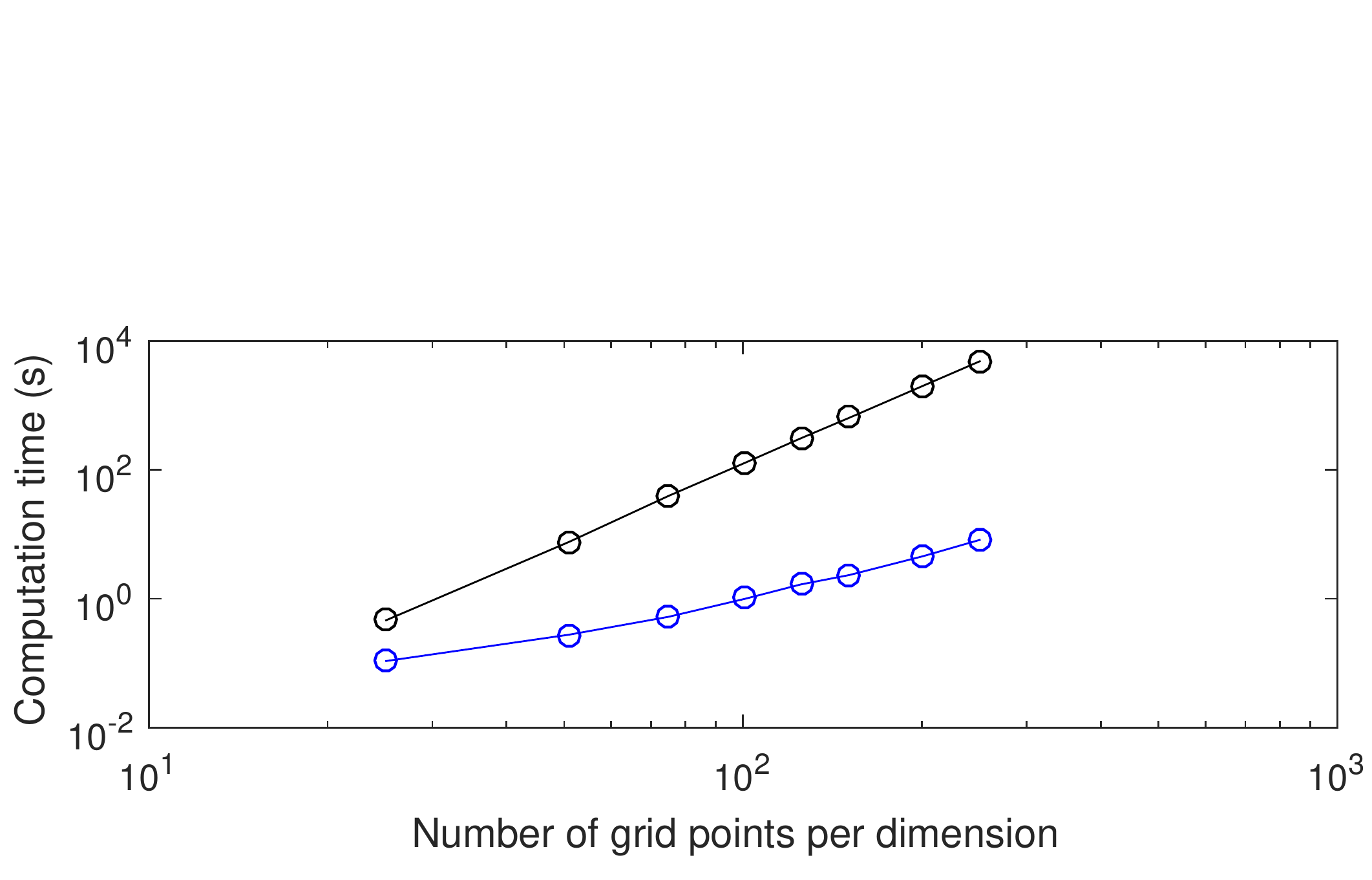}
	\caption{Computation times of the two methods in log scale for the Dubins Car. The time of the direct computation in 3D increases rapidly with the number of grid points per dimension. In contrast, computation times in 2D with decomposition are negligible in comparison.}
	\label{fig:dubins_time}
\end{figure}

The computation benefits of using our decomposition method can be seen from Fig. \ref{fig:dubins_time}. The plot shows, in log-log scale, the computation time in seconds versus the number of grid points per dimension in the numerical computation. One can see that the direct computation of the BRS in 3D becomes very time-consuming as the number of grid points per dimension is increased, while the computation via decomposition hardly takes any time in comparison. Directly computing the BRS with 251 grid points per dimension in 3D took approximately 80 minutes, while computing the BRS via decomposition in 2D only took approximately 30 seconds! The computations were timed on a computer with an Intel Core i7-2600K processor and 16GB of random-access memory.

\begin{figure}
	\centering
	\includegraphics[width=0.9\columnwidth]{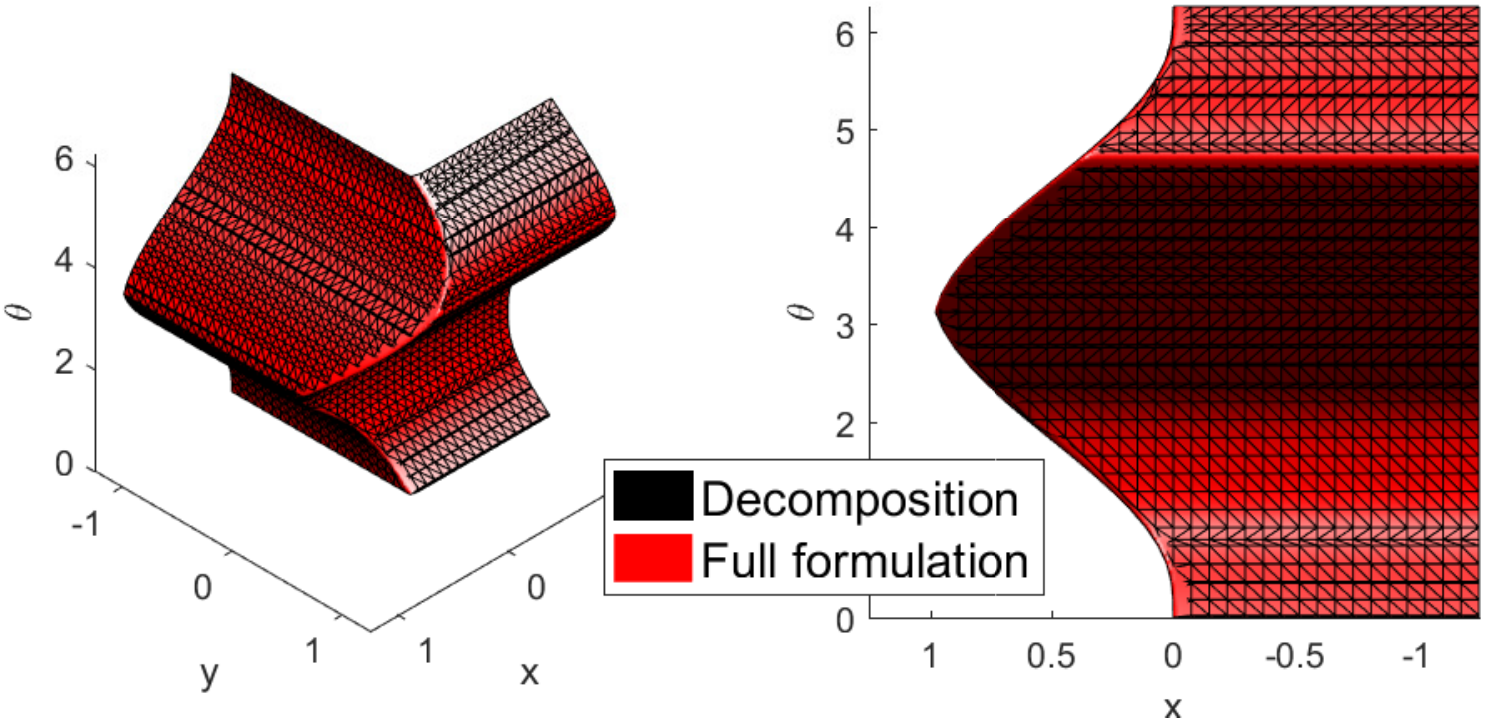}
	\caption{Comparison of the $\maxbrs(t)$ computed using our decomposition method and the full formulation. The computation results are indistinguishable. Note that the surface shows the boundary of the set; the set itself is on the ``near'' side of the left subplot, and the left side of the right subplot.}
	\label{fig:dubins_un_set}
\end{figure}

Fig. \ref{fig:dubins_un_set} illustrates Theorem \ref{thm:sc_reach_un}. We chose the target set to be $\targetset = \{(\pos_x, \pos_y, \theta): \pos_x \le 0.5 \vee \pos_y \le 0.5\}$, and computed the BRS $\maxbrs(t), t = -0.5$ via decomposition. No additional approximation error is incurred in the reconstruction process. The target set can be written as $\targetset = \bp(\sctarget{1}) \cup \bp(\sctarget{1})$ where $\sctarget{1} = \{(\pos_x, \theta): \pos_x \le 0.5\}, \sctarget{2} = \{(\pos_y, \theta ): \pos_y \le 0.5\}$.


\section{SCSs with Decoupled Control\label{sec:decoupled}}
In this section, we consider a special case of \eqref{eq:scdyn} in which the subsystem controls do not have any shared components. The results from \ref{sec:sc} still hold, and in addition we can state the results in Propositions \ref{prop:decoupled_reach_in} and \ref{prop:decoupled_reach_un}. The special case of \eqref{eq:scdyn} is as follows:

\begin{equation}
\label{eq:ddyn}
\begin{aligned}
\dot\spart_1 &= \fdyn_1(\spart_1, \spart_c, \cpart_1) \\
\dot\spart_2 &= \fdyn_2(\spart_2, \spart_c, \cpart_2) \\
\dot\spart_c &= \fdyn_c(\spart_c) \\
\end{aligned}
\end{equation}

%

\noindent where the subsystem controls do not have any shared components, so that we have $\ctrl = (\cpart_1, \cpart_2)$. Furthermore, we can define the trajectory of $\spart_c$ as $\eta(s; \spart_c, t))$, which satisfies

\begin{equation}
\begin{aligned}
\frac{d}{ds}\eta_c(s; \spart_c, t) &= \fdyn_c(\eta_c(s; \spart_c, t)) \\
\eta_c(t; \spart_c, t)) &= \spart_c
\end{aligned}
\end{equation}

Note that since the trajectory $\eta_c(s; \spart_c, t)$ does not depend on the control, we can treat $\eta_c(s; \spart_c, t)$ as a constant when given $\spart_c$ and $s$. Therefore, given $\spart_c$, the other \textit{state partitions} also become self-contained, with dynamics

\begin{equation}
\dot\spart_i = \fdyn_i(\spart_i, \spart_c, \cpart_i) = \fdyn_i(\spart_i, \cpart_i; \eta_c(s; \spart_c, t))
\end{equation}

\noindent and with trajectories $\eta_i(s; \spart_i, \spart_c, t, \cpart_i(\cdot))$ satisfying

\begin{equation}
\begin{aligned}
&\frac{d}{ds}\eta_i(s; \spart_i, \spart_c, t, \cpart_i(\cdot)) \\
&\qquad = \fdyn_i(\eta_i(s; \spart_i, \spart_c, t, \cpart_i(\cdot)), \cpart_i(s); \eta_c(s; \spart_c, t)) \\
&\qquad \eta_i(t; \spart_i, \spart_c, t, \cpart_i(\cdot)) = \spart_i
\end{aligned}
\end{equation}

Therefore, the \textit{subsystem trajectories} can be written as

\begin{equation}
\label{eq:sctraj_decoupled}
\sctraj_i(s; \scstate_i, t, \cpart_i(\cdot)) = \big(\eta_i(s; \spart_i, \spart_c, t, \cpart_i(\cdot)), \eta_c(s; \spart_c, t)\big)
\end{equation}

\begin{prop}
\label{prop:decoupled_reach_in}
Suppose that the full system in \eqref{eq:fdyn} can be decomposed into the form of \eqref{eq:ddyn}. Then,

\begin{equation}
\label{eq:targetset_in_decoupled}
\begin{aligned}
&\targetset = \bp(\sctarget{1}) \cap \bp(\sctarget{2}) \\
&\Rightarrow \maxbrs(t) = \bp(\maxbrs_1(t)) \cap \bp(\maxbrs_2(t))
\end{aligned}
\end{equation}
\end{prop}

\begin{prop}
\label{prop:decoupled_reach_un}
Suppose that the full system in \eqref{eq:fdyn} can be decomposed into the form of \eqref{eq:ddyn}. Then,

\begin{equation}
\label{eq:targetset_un_decoupled}
\begin{aligned}
&\targetset = \bp(\sctarget{1}) \cup \bp(\sctarget{2}) \\
&\Rightarrow \minbrs(t) = \bp(\minbrs_1(t)) \cup \bp(\minbrs_2(t))
\end{aligned}
\end{equation}
\end{prop}

\begin{rem}
    Systems with fully decoupled subsystems in the form of $\scstate_1 = \spart_1, \scstate_2 = \spart_2$ are a special case of \eqref{eq:ddyn}. A numerical example illustrating this case will be presented in the 10D Near-Hover Quadrotor example in Section \ref{sec:hd10}.
\end{rem}

%
%
%
%
%

\subsection{Proof of Proposition \ref{prop:decoupled_reach_in}}
We will prove the following equivalent statement:
\begin{equation}
\label{eq:proof_left}
\bar\state \in \maxbrs(t) \Leftrightarrow \bar\state \in \bp(\maxbrs_1(t)) \cap \bp(\maxbrs_2(t))
\end{equation}

By Definition \ref{defn:rset_goal} (maximal BRS), we have 

\begin{equation}
\label{eq:decoupled_in}
\bar\state \in \maxbrs(t)\Leftrightarrow \exists \ctrl(\cdot), \traj(0; \bar\state, t, \ctrl(\cdot)) \in \targetset \\
\end{equation}

Consider the relationship between the full system trajectory and subsystem trajectory in \eqref{eq:proj_traj}. Define
\begin{equation*}
\begin{aligned}
&\bar \scstate_i = (\bar\spart_i, \bar\spart_c) = \proj_i(\bar\state), \text{ and} \\
&\sctraj_i(s; \bar \scstate_i, t, \cpart_i(\cdot)) = \proj_i(\traj(0; \bar\state, t, \ctrl(\cdot)))
\end{aligned}
\end{equation*}

By \eqref{eq:sctraj_decoupled} we can write

\begin{equation*}
\big(\eta_i(s; \bar\spart_i, \bar\spart_c, t, \cpart_i(\cdot)), \eta_c(s; \bar\spart_c, t)\big) = \proj_i(\traj(0; \bar\state, t, \ctrl(\cdot)))
\end{equation*}

Since $\targetset = \bp(\sctarget{1}) \cap \bp(\sctarget{2})$, \eqref{eq:decoupled_in} is equivalent to

\begin{equation}
\label{eq:decoupled_ctrl_before_distr}
\begin{array}{rc}
&
\begin{aligned}
&\exists \left(\cpart_1(\cdot), \cpart_2(\cdot)\right), \\
&\quad\big(\eta_1(s; \bar\spart_1, \bar\spart_c, t, \cpart_1(\cdot)), \eta_c(s; \bar\spart_c, t)\big) \in \sctarget{1} \wedge\\
&\quad\big(\eta_2(s; \bar\spart_2, \bar\spart_c, t, \cpart_2(\cdot)), \eta_c(s; \bar\spart_c, t)\big) \in \sctarget{2}
\end{aligned} \\
&\text{(by Corollary \ref{cor:proj_in})}
\end{array}
\end{equation}

\begin{equation}
\label{eq:decoupled_ctrl_after_distr}
\begin{array}{rc}
\Leftrightarrow
&
\begin{aligned}
&\exists \cpart_1(\cdot), \big(\eta_1(s; \bar\spart_1, \bar\spart_c, t, \cpart_1(\cdot)), \eta_c(s; \bar\spart_c, t)\big) \in \sctarget{1} \wedge\\
&\exists \cpart_2(\cdot), \big(\eta_2(s; \bar\spart_2, \bar\spart_c, t, \cpart_2(\cdot)), \eta_c(s; \bar\spart_c, t)\big) \in \sctarget{2}
\end{aligned} \\
&\text{(since subsystem controls do not share components)}
\end{array}
\end{equation}

\begin{equation*}
\begin{array}{rc}
\Leftrightarrow
& 
\scstate_1 \in \maxbrs_1(t) \wedge \scstate_2 \in \maxbrs_2(t) \\
&\text{(by definition of subsystem BRS in \eqref{eq:sc_rset})}
\end{array}
\end{equation*}

\begin{equation*}
\begin{array}{rc}
\Leftrightarrow
& 
\bar\state \in \bp(\maxbrs_1(t)) \cap \bp(\maxbrs_2(t)) \\
&\text{(by Corollary \ref{cor:proj_in})}
\end{array}
\end{equation*}
\hfill\IEEEQEDhere

\subsection{Proof of Proposition \ref{prop:decoupled_reach_un}}
This proof follows the same arguments as \ref{prop:decoupled_reach_in}, but with ``$\maxbrs$'', ``$\cap$'',  ``$\exists$'' replaced with ``$\minbrs$'', ``$\cup$'', ``$\forall$'', respectively. \hfill \IEEEQEDhere

\begin{rem}
When the subsystem controls have shared components, the control chosen by each subsystem may not agree with the other. This is the intuition behind why the results of Propositions \ref{prop:decoupled_reach_in} and \ref{prop:decoupled_reach_un} only hold true when there are no shared components in the subsystem controls. Note that the theorems hold despite the state coupling between the subsystems.
\end{rem}
The results from this section are summarized in Table \ref{table:Decoupled_BRS}.

\begin{table}[h]
	\centering
	\caption{BRS Results from Section \ref{sec:decoupled}}
	\label{table:Decoupled_BRS}
	\SHnote{\begin{tabular}{|l|l|l|}
		\hline
		\textbf{Shared Controls}       & \multicolumn{2}{c|}{\textbf{No}}                                                                                                                 \\ \hline
		\textbf{Shared Disturbance}    & \multicolumn{2}{c|}{\textbf{No}}                                                                                                                 \\ \hline
		\textbf{Target}                & \textbf{Intersection}                                                       & \textbf{Union}                                                              \\ \hline
		Recover Max. BRS?     & Yes, exact                                                         & Yes, exact                                                         \\ \hline
		Recover Min. BRS?     & Yes, exact                                                         & Yes, exact                                                         \\ \hline
		Equation(s) & \begin{tabular}[c]{@{}l@{}}Prop \ref{prop:decoupled_reach_in}, (\ref{eq:targetset_in_decoupled})\\ Thm \ref{thm:sc_reach_in}, (\ref{eq:targetset_in})\end{tabular} & \begin{tabular}[c]{@{}l@{}}Thm \ref{thm:sc_reach_un}, (\ref{eq:targetset_un})\\ Prop \ref{prop:decoupled_reach_un}, (\ref{eq:targetset_un_decoupled})\end{tabular} \\ \hline
	\end{tabular}}
\end{table} 

\section{Decomposition of Reachable Tubes \label{sec:set_to_tube}}
Sometimes, BRTs are desired; for example, in safety analysis, the computation of the BRT $\minbrt(t)$ in Definition \ref{defn:rset_avoid} is quite important, since if the target set $\targetset$ represents an unsafe set of states, then $\minbrt(t)$ contains all states that would lead to some unsafe state at \textit{some time} within a duration of length $|t|$. 

We now first discuss a special case where the full system BRT can be directly reconstructed from subsystem BRTs in Section \ref{subsec:tube2tube}, and then present a general method in which a BRT can be obtained via the union of BRSs in Section \ref{subsec:union_BRT}.

\subsection{Full System BRTs From Subsystem BRTs \label{subsec:tube2tube}}
Intuitively, it may seem like the results related to BRSs outlined in Sections \ref{sec:sc} and \ref{sec:decoupled} trivially carry over to BRTs, and the relationship between BRSs and BRTs are relatively simple; however, this is only partially true. The results related to BRSs presented so far in this paper only easily carry over for BRTs if $\targetset = \bp(\sctarget{1}) \cup \bp(\sctarget{2})$. This is formally stated in the following Proposition:

\begin{prop}
\label{prop:BRT_Union}
Suppose \eqref{eq:target_union} holds, that is, 
\begin{equation*}
\targetset = \bp(\sctarget{1}) \cup \bp(\sctarget{2})
\end{equation*}

Then, the full-dimensional BRT can be reconstructed from the lower-dimensional BRTs without incurring additional approximation errors. For systems with SCSs as in \eqref{eq:scdyn}, we have

\begin{equation}
\label{eq:maxBRT_union}
\maxbrt(t) = \bp(\maxbrt_1(t)) \cup \bp(\maxbrt_2(t))
\end{equation}

This can be proven by starting the proof of Theorem \ref{thm:sc_reach_un} with Definition \ref{defn:rtube_goal} for the BRT $\maxbrt(t)$ instead of Definition \ref{defn:rset_goal} for the BRS $\maxbrs(t)$.

For systems with subsystem controls that do not share any components, we \textit{in addition} have
\begin{equation}
\label{eq:minBRT_union}
\minbrt(t) = \bp(\minbrt_1(t)) \cup \bp(\minbrt_2(t))
\end{equation}

This can be proven by starting the proof of Proposition \ref{prop:decoupled_reach_un} from Definition \ref{defn:rtube_avoid} for the BRT $\minbrt(t)$ instead of Definition \ref{defn:rset_avoid} for the BRS $\minbrs(t)$. Note that (\ref{eq:minBRT_union}) does not necessarily hold for systems with shared controls.
\end{prop}

\subsection{BRTs From Union of BRSs \label{subsec:union_BRT}}
If $\targetset = \bp(\sctarget{1}) \cap \bp(\sctarget{2})$, the BRT cannot be directly reconstructed from lower-dimensional BRTs because when computing with BRTs, we lose information about the exact time that a trajectory enters a set. Instead, we provide a general method of obtaining the BRT that can \textit{also} be used in the case of Section \ref{subsec:tube2tube}: We first compute the BRSs, and then take their union to obtain the BRT. For this case, we show that $\maxbrt(t) = \bigcup_{s\in[t,0]} \maxbrs(s)$, and $\minbrt(t) = \bigcup_{s\in[t,0]} \minbrs(s)$ when $\minbrs(s) \neq \emptyset ~ \forall s \in [t,0]$. These results related to the indirect reconstruction of BRTs are given in Proposition \ref{prop:sets2tube_goal} and Theorem \ref{thm:sets2tube_avoid}

\begin{prop}
\label{prop:sets2tube_goal}
\begin{equation}
\label{eq:maxBRT_intersection}
\bigcup_{s \in [t, 0]} \maxbrs(s) = \maxbrt(t)
\end{equation}
\end{prop}

\begin{thm}
\label{thm:sets2tube_avoid}
\begin{equation}
\bigcup_{s \in [t, 0]} \minbrs(s) \subseteq \minbrt(t)
\end{equation}

In addition, if $\forall s \in [t, 0], \minbrs(s) \neq \emptyset$, then
\begin{equation}
\label{eq:minBRT_intersection}
\bigcup_{s \in [t, 0]} \minbrs(s) = \minbrt(t).
\end{equation}
\end{thm}

Propositions \ref{prop:sets2tube_goal} and the first part of Theorem \ref{thm:sets2tube_avoid} are known \cite{Mitchell07b}, but we present them in our paper in greater detail for clarity and completeness. The second part of Theorem \ref{thm:sets2tube_avoid} is the main new result related to obtaining the BRT from BRSs. 

\begin{rem}
The reason the theorems in Sections \ref{sec:sc}  and \ref{sec:decoupled} trivially carry over when $\targetset = \bp(\sctarget{1}) \cup \bp(\sctarget{2})$ is that in this case, any subsystem trajectory that reaches the corresponding subsystem target set implies that the full system trajectory reaches the full system target set. 

In contrast, in the case $\targetset = \bp(\sctarget{1}) \cap \bp(\sctarget{2})$, \textit{both} subsystem trajectories must be in the corresponding subsystem target sets \textit{at the same time}. Mathematically, recall the definitions of subsystem BRTs in \eqref{eq:sc_rtube}:

\begin{equation*}
\begin{aligned}
\minbrt_i(t) = \{&\scstate_i: \forall \ctrl(\cdot), \exists s \in [t, 0], \sctraj_i(s; \scstate_i, t, \ctrl(\cdot)) \in \sctarget{i} \} \\
\maxbrt_i(t) = \{&\scstate_i: \exists \ctrl(\cdot), \exists s \in [t, 0],\sctraj_i(s; \scstate_i, t, \ctrl(\cdot)) \in \sctarget{i} \} \\
\end{aligned}
\end{equation*}

The set of ``$s$'' during which each subsystem trajectory is in $\sctarget{i}$ may not overlap for the different subsystems. In this case, we can still first compute the BRSs in lower-dimensional subspaces, and then convert the BRSs to the BRT using Propositions \ref{prop:sets2tube_goal} and Theorem \ref{thm:sets2tube_avoid}.
\end{rem}

\subsection{Proof of Proposition \ref{prop:sets2tube_goal}}
\label{sec:sets2tube_goal}
We start with Definition \ref{defn:rset_goal} (maximal BRS):
\begin{equation*}
\maxbrs(t) = \{\state: \exists \ctrl(\cdot) \in \cfset, \traj(0; \state, t, \ctrl(\cdot)) \in \targetset \}
\end{equation*}

If some state $\state$ is in the union $\bigcup_{s \in [t, 0]} \maxbrs(s)$, then there is some $s \in [t,0]$ such that $\state\in\maxbrs(s)$. Therefore, the union can be written as

\begin{equation}
\label{eq:union_goal_expression}
\bigcup_{s \in [t, 0]}\maxbrs(s) = \{\state: \exists s \in [t, 0], \exists \ctrl(\cdot), \traj(0; \state, s, \ctrl(\cdot)) \in \targetset \}
\end{equation}

Suppose $\state \in \bigcup_{s \in [t, 0]}\maxbrs(s)$, then equivalently

\begin{equation}
\label{eq:union_avoid_state}
\exists s \in [t, 0], \exists \ctrl(\cdot) \in \cfset, \traj(0; \state, s, \ctrl(\cdot)) \in \targetset
\end{equation}

Using \eqref{eq:time-invariant}, the time-invariance of the system, we can shift the trajectory time arguments by $t - s$ to get

\begin{equation}
\exists s \in [t, 0], \exists \ctrl(\cdot) \in \cfset, \traj(t-s; \state, t, \ctrl(\cdot)) \in \targetset
\end{equation}

Since $s \in [t, 0] \Leftrightarrow t-s \in [t, 0]$, we can equivalently write

\begin{equation}
\exists s \in [t, 0], \exists \ctrl(\cdot) \in \cfset, \traj(s; \state, t, \ctrl(\cdot)) \in \targetset
\end{equation}

We can swap the expressions $\exists s \in [t, 0]$ and $\exists \ctrl(\cdot) \in \cfset$ without changing meaning since both quantifiers are the same:

\begin{equation}
\label{eq:goal_after_swap}
\exists \ctrl(\cdot) \in \cfset, \exists s \in [t, 0], \traj(s; \state, t, \ctrl(\cdot)) \in \targetset
\end{equation}

\noindent which is equivalent to $\state \in \maxbrt(t)$ by Definition \ref{defn:rtube_goal} (maximal BRT). \hfill \IEEEQEDhere

\subsection{Proof of Theorem \ref{thm:sets2tube_avoid}}
\label{sec:sets2tube_avoid}

We first establish $\bigcup_{s \in [t, 0]} \minbrs(s) \subseteq \minbrt(t)$. Consider Definition \ref{defn:rset_avoid} (minimal BRS):
\begin{equation*}
  \minbrs(t) = \{\state: \forall \ctrl(\cdot) \in \cfset, \traj(0; \state, t, \ctrl(\cdot)) \in \targetset \}
\end{equation*}

If some state $\state$ is in the union $\bigcup_{s \in [t, 0]} \minbrs(s)$, then $\exists s \in [t,0]$ such that $\state\in\minbrs(s)$. Thus, the union can be written as

\begin{equation}
  \label{eq:union_avoid_expression}
  \bigcup_{s \in [t, 0]}\minbrs(s) = \{\state: \exists s \in [t, 0], \forall \ctrl(\cdot), \traj(0; \state, s, \ctrl(\cdot)) \in \targetset \}
\end{equation}

Suppose $\state \in \bigcup_{s \in [t, 0]}\minbrs(s)$, then

\begin{equation}
  \label{eq:union_avoid_state}
  \exists s \in [t, 0], \forall \ctrl(\cdot) \in \cfset, \traj(0; \state, s, \ctrl(\cdot)) \in \targetset
\end{equation}

Using \eqref{eq:time-invariant}, the time-invariance of the system, we can shift the trajectory time arguments by $t - s$ to get

\begin{equation}
  \exists s \in [t, 0], \forall \ctrl(\cdot) \in \cfset, \traj(t-s; \state, t, \ctrl(\cdot)) \in \targetset
\end{equation}

Since $s \in [t, 0] \Leftrightarrow t-s \in [t, 0]$, we can equivalently write

\begin{equation}
  \exists s \in [t, 0], \forall \ctrl(\cdot) \in \cfset, \traj(s; \state, t, \ctrl(\cdot)) \in \targetset
\end{equation}

Let such an $s \in [t,0]$ be denoted $\bar s$, then

\begin{equation}
  \begin{aligned}
    &\forall \ctrl(\cdot) \in \cfset, \traj(\bar s; \state, t, \ctrl(\cdot)) \in \targetset\\
    &\Rightarrow\forall \ctrl(\cdot) \in \cfset, \exists s \in [t, 0], \traj(s; \state, t, \ctrl(\cdot)) \in \targetset
  \end{aligned}
\end{equation}

By Definition \ref{defn:rtube_avoid}, we have $\state \in \minbrt(t)$.

Next, given $\forall s \in [t, 0], \minbrs(s) \neq \emptyset$, we show $\bigcup_{s \in [t, 0]} \minbrs(s) \supseteq \minbrt(t)$. Equivalently, we show

\begin{equation}
\state \notin \bigcup_{s \in [t, 0]} \minbrs(s) \Rightarrow \state \notin \minbrt(t)
\end{equation}

%

First, observe that by the definition of minimal BRS, we have that if any state $\bar\state \in \minbrs(t)$, then

\begin{equation}
\forall s \in [t,0], \forall \ctrl(\cdot) \in \cfset, \traj(s; \bar\state, t, \ctrl(\cdot)) \in \minbrs(s)
\end{equation}

\noindent since otherwise, we would have for some $\bar s \in [t,0]$,

\begin{equation}
\begin{aligned}
&\exists \ctrl(\cdot) \in \cfset, \traj(\bar s; \bar\state, t, \ctrl(\cdot)) \notin \minbrs(\bar s) \\
&\Rightarrow \exists \ctrl(\cdot) \in \cfset, \traj(0; \traj(\bar s; \bar\state, t, \ctrl(\cdot)), \bar s, \ctrl(\cdot)) \notin \targetset \\
&\Leftrightarrow \exists \ctrl(\cdot) \in \cfset, \traj(0; \bar\state, t, \ctrl(\cdot)) \notin \targetset
\end{aligned}
\end{equation}

\noindent which contradicts $\bar\state \in \minbrs(t)$.
%
%
%
%
%
%
%
%
%
%

Given $\state \notin \minbrs(t)$, there exists some control $\bar \ctrl(\cdot)$ such that $\traj(0; \state, t, \bar\ctrl(\cdot)) \notin \targetset = \minbrs(0)$. Moreover, we must have $\forall s \in [t,0], \traj(s; \state, t, \bar\ctrl(\cdot)) \notin \targetset$, since otherwise, we would have $\exists \hat s$ such that

\begin{equation}
\begin{aligned}
&\traj(\hat s; \state, t, \bar\ctrl(\cdot)) \in \targetset =\minbrs(0) \\
&\Rightarrow\state = \traj(t; \state, t, \bar\ctrl(\cdot)) \in \minbrs(t - \hat s)
\end{aligned}
\end{equation}

\noindent which contradicts $\state \notin \bigcup_{s \in (t,0)} \minbrs(s)$.

Using time-invariance of the system dynamics, we have $\forall s \in [t,0], \traj(0; \state, t-s, \bar\ctrl(\cdot)) \notin \targetset$, which is equivalent to $\forall s \in [t,0], \traj(0; \state, s, \bar\ctrl(\cdot)) \notin \targetset$. Therefore, $\exists \ctrl(\cdot) \in \cfset, \forall s \in [t,0], \traj(0; \state, s, \bar\ctrl(\cdot)) \notin \targetset \Leftrightarrow \state \notin \minbrt(t)$. \hfill \IEEEQEDhere
%
%
%
%
%
%
%

\begin{rem}
When $\exists s \in [t, 0], \minbrs(s) = \emptyset$, it is currently not known whether the union of the BRSs $\minbrs(s)$ will be equal to the BRT $\minbrt(t)$ or a proper subset of the BRT $\minbrt(t)$. Both are possibilities. Finding a weaker condition under which the union of BRSs equals to the BRT is an important future direction that we plan to investigate.
\end{rem}

\begin{rem}
Note that Proposition \ref{prop:sets2tube_goal} and Theorem \ref{thm:sets2tube_avoid} also hold for decoupled control.
\end{rem}

\begin{table}[h]
	\centering
	\caption{BRT Results for Reconstruction from Tubes}
	\label{table:BRT_Tubes}
	\SHnote{{\scriptsize
	\begin{tabular}{|l|l|l|l|l|}
		\hline
		\textbf{Shared Controls}    & \multicolumn{2}{c|}{\textbf{Yes}}  & \multicolumn{2}{c|}{\textbf{No}}   \\ \hline
		\textbf{Shared Dstb.} & \multicolumn{2}{c|}{\textbf{No}}   & \multicolumn{2}{c|}{\textbf{No}}   \\ \hline
		\textbf{Target}             & \textbf{Intersection} & \textbf{Union}      & \textbf{Intersection} & \textbf{Union}      \\ \hline
		\begin{tabular}[c]{@{}l@{}}Recover \\ Max. BRT?  \end{tabular} & No           & Yes, exact & No           & Yes, exact \\ \hline
		\begin{tabular}[c]{@{}l@{}}Recover \\ Min. BRT?  \end{tabular}  & No           & No         & No           & Yes, exact \\ \hline
		Equation(s)        &     \begin{tabular}[c]{@{}l@{}}N/A, see\\ Table \ref{table:BRT_Sets}  \end{tabular}      &  \begin{tabular}[c]{@{}l@{}}Prop \ref{prop:BRT_Union},\\ (\ref{eq:maxBRT_union})  \end{tabular}            &        \begin{tabular}[c]{@{}l@{}}N/A, see\\ Table \ref{table:BRT_Sets}  \end{tabular}      & \begin{tabular}[c]{@{}l@{}}Prop \ref{prop:BRT_Union}, (\ref{eq:maxBRT_union})\\ Prop \ref{prop:BRT_Union}, (\ref{eq:minBRT_union})\end{tabular}        \\ \hline
	\end{tabular}
}}
\end{table}

\begin{table}[h]
	\centering
	\caption{BRT Results for Reconstruction from Sets}
	\label{table:BRT_Sets}
	\SHnote{{\scriptsize
	\begin{tabular}{|l|c|l|c|l|}
		\hline
		Shared Controls    & \multicolumn{2}{c|}{Yes}                  & \multicolumn{2}{c|}{No}                   \\ \hline
		Shared Disturbance & \multicolumn{2}{c|}{No}                   & \multicolumn{2}{c|}{No}                   \\ \hline
		Target             & \multicolumn{1}{l|}{Intersection} & Union & \multicolumn{1}{l|}{Intersection} & Union \\ \hline
		Recover Max. BRT?  & \multicolumn{4}{c|}{Yes, exact}                                                       \\ \hline
		Recover Min. BRT?  & \multicolumn{4}{c|}{Yes, exact*}                                                      \\ \hline
		Equation(s)        & \multicolumn{4}{c|}{\begin{tabular}[c]{@{}l@{}}Prop \ref{prop:sets2tube_goal}, (\ref{eq:maxBRT_intersection})\\ Thm \ref{thm:sets2tube_avoid}, (\ref{eq:minBRT_intersection})\end{tabular}}                                                                \\ \hline
	\end{tabular}
}}
\end{table}

\subsection{Numerical Results}

\begin{figure}[H]
	\centering
	\includegraphics[width=0.9\columnwidth]{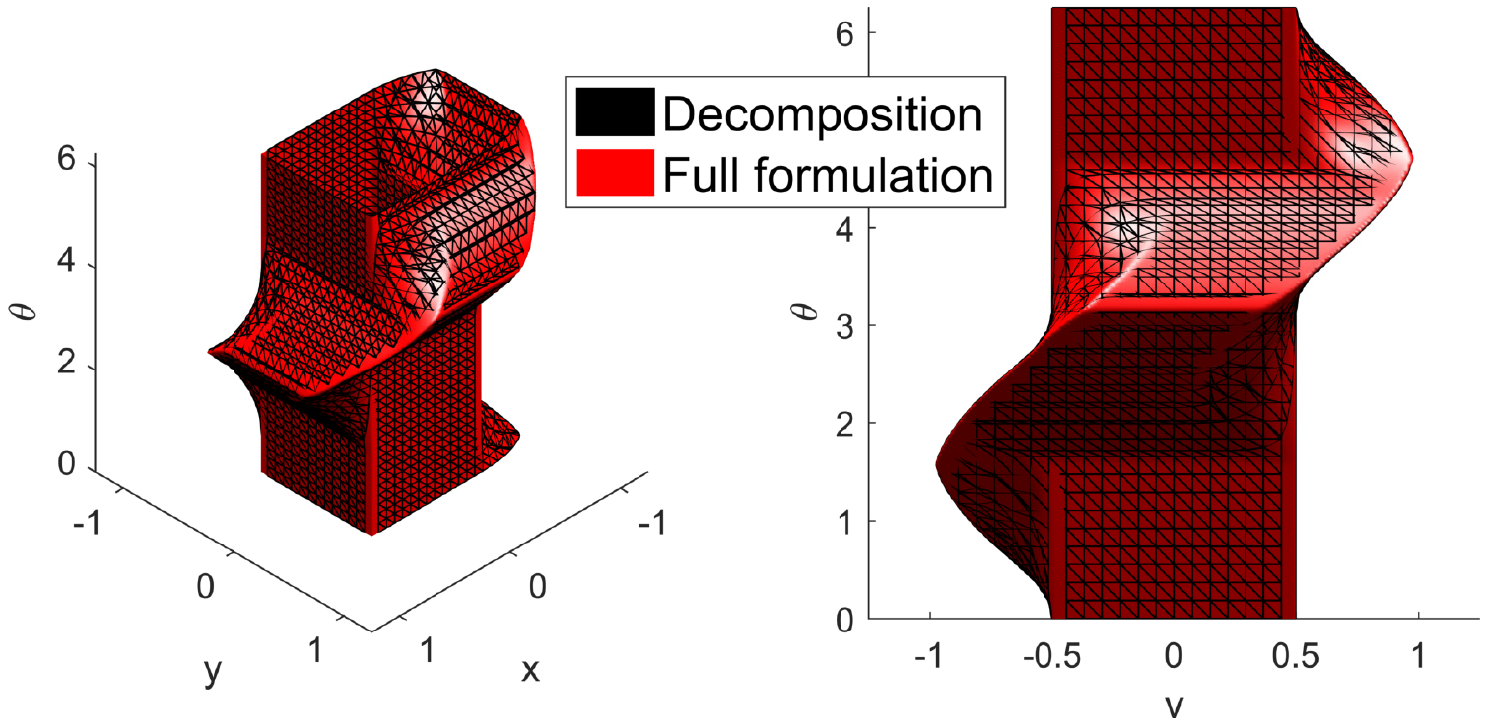}
	\caption{The BRT computed directly in 3D (red surface) and computed via decomposition in 2D (black mesh). Using our decomposition techinque, we first compute the BRSs $\minbrs(s), s\in[-0.5,0]$, and then obtained the BRT by taking their union.}
	\label{fig:dubins_in_tube}
\end{figure}

We now revisit the Dubins Car, whose full system and subsystem dynamics are given in \eqref{eq:dubins_car} and \eqref{eq:dubins_car_decomp} respectively. Using the target set $\targetset$ given in \eqref{eq:dubins_target} and writing $\targetset$ in the form of \eqref{eq:dubins_target_decomp}, we computed the BRT $\minbrt(t), t=-0.5$ by first computing $\minbrs(s), s\in[-0.5, 0]$, and then taking their union.

Fig. \ref{fig:dubins_in_tube} shows the BRT $\minbrt(t), t = -0.5$ computed directly in 3D and via decomposition. Since $\minbrs(s) \neq \emptyset ~ \forall s\in[-0.5,0]$, the reconstruction does not incur any additional approximation errors. 

\section{High-Dimensional Numerical Results \label{sec:hd}}
In this section, we show numerical results for the 6D Acrobatic Quadrotor and the 10D Near-Hover Quadrotor, two systems whose exact BRSs and BRTs were intractable to compute with previous methods to the best of our knowledge.

\subsection{The 6D Acrobatic Quadrotor}
In \cite{Gillula11}, a 6D quadrotor model used to perform backflips was simplified into a series of smaller models linked together in a hybrid system. The Quadrotor has state $\state = (\pos_x, \vel_x, \pos_y, \vel_y, \phi, \omega)$, and dynamics

\begin{equation}
\label{eq:quad6D}
\begin{array}{c}
	\left[
	\begin{array}{c}
	\dot\pos_x\\
	\dot\vel_x\\
	\dot\pos_y\\
	\dot\vel_y\\
	\dot\phi\\
	\dot\omega
	\end{array}
	\right]
	=
	\left[
	\begin{array}{c}
	\vel_x\\
	-\frac{1}{m}C^v_D\vel_x -\frac{T_1}{m}\sin\phi  -\frac{T_2}{m}\sin\phi \\
	\vel_y\\
	-\frac{1}{m}\left( mg+C^v_D\vel_y\right) + \frac{T_1}{m}\cos\phi + \frac{T_2}{m}\cos\phi\\
	\omega\\
	-\frac{1}{I_{yy}}C^\phi_D\omega -\frac{l}{I_{yy}}T_1 + \frac{l}{I_{yy}}T_2
	\end{array}
	\right]\\
	\end{array}
\end{equation}
where $x$, $y$, and $\phi$ represent the quadrotor's horizontal, vertical, and rotational positions, respectively. Their derivatives represent the velocity with respect to each state. The control inputs $T_1$ and $T_2$ represent the thrust exerted on either end of the quadrotor, and the constant system parameters are $m$ for mass, $C_D^v$ for translational drag, $C_D^\phi$ for rotational drag, $g$ for acceleration due to gravity, $l$ for the length from the quadrotor's center to an edge, and $I_{yy}$ for moment of inertia. 

We decompose the system into the following subsystems:
\begin{equation}
\label{eq:quad6Dsc}
\scstate_1 = (\pos_x, \vel_x, \phi, \omega), \qquad \scstate_2 = (\pos_y, \vel_y, \phi, \omega)
\end{equation}

\begin{figure}
	\centering
	\includegraphics[width=\columnwidth]{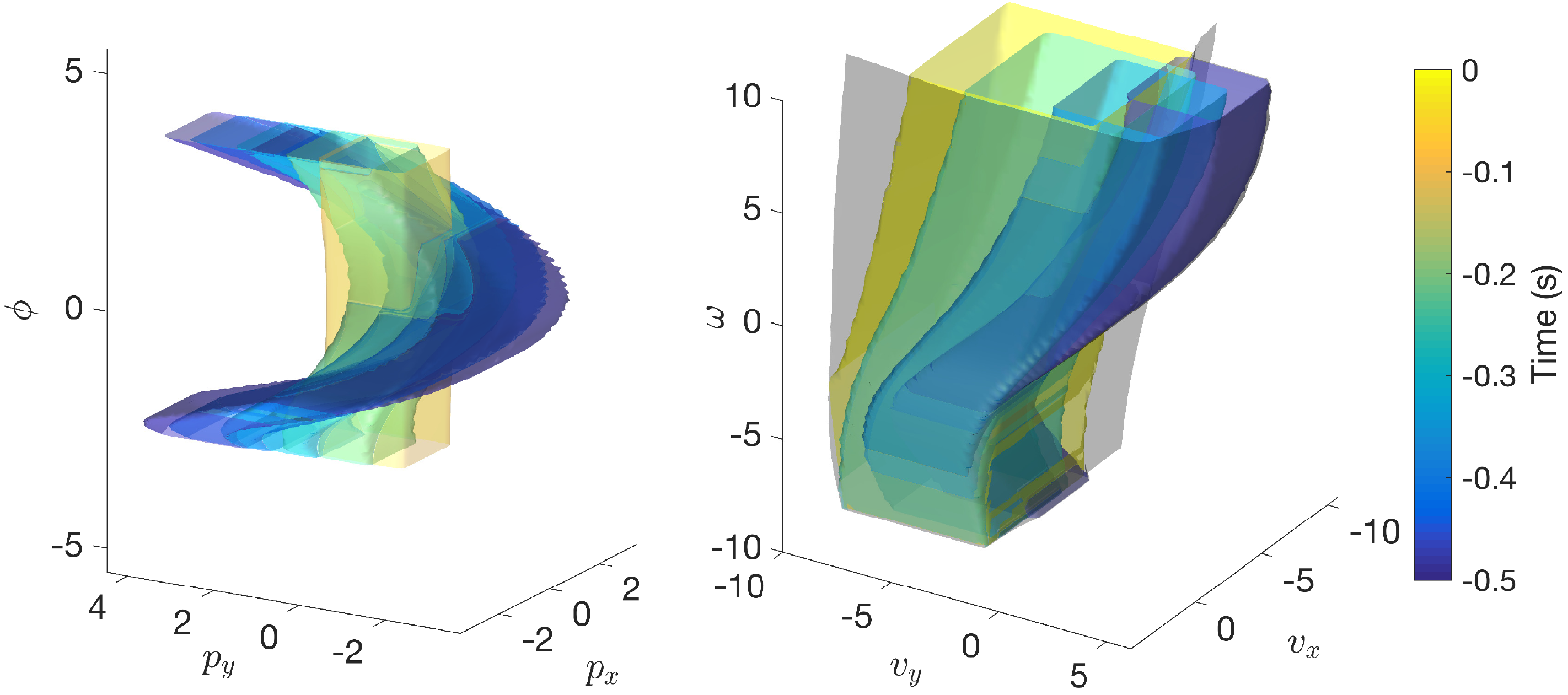}
	\caption{Left: 3D positional slices of the reconstructed 6D BRSs at $\vel_x=\vel_y=1$, $\omega=0$ at different points in time. The BRT cannot be seen in this image because it encompasses the entire union of BRSs. Right: 3D velocity slices of the reconstructed 6D BRSs at $x,y=1.5$, $\phi = 1.5$ at different points in time. The BRT can be seen as the transparent gray surface that encompasses the sets.}
	\label{fig:Quad6D_Combo}
\end{figure}

For this example we will compute $\minbrs(t)$ and $\minbrt(t)$, which describe the set of initial conditions from which the system may enter the target set despite the best possible control to avoid the target. We define the target set as a square of length 2 centered at $(\pos_x,\pos_y)=(0,0)$ described by $\targetset = \{(\pos_x,\vel_x, \pos_y, \vel_y, \phi, \omega): |\pos_x|, |\pos_y| \le 1\}$. This can be interpreted as a positional box centered at the origin that must be avoided for all angles and velocities. From the target set, we define $\fc(\state)$ such that $\fc(\state)\le 0 \Leftrightarrow \state\in\targetset$. This target set is then decomposed as follows:
\begin{equation*}
\begin{aligned}
\sctarget{1} &= \{(\pos_x, \vel_x,\phi, \omega): |\pos_x|\le 1\} \\
\sctarget{2} &= \{(\pos_y, \vel_y, \phi, \omega): |\pos_y|\le 1\}
\end{aligned}
\end{equation*}
The BRS of each 4D subsystem is computed and then recombined into the 6D BRS. To visually depict the 6D BRS, 3D slices of the BRS along the positional and velocity axes were computed. The left image in Fig. \ref{fig:Quad6D_Combo} shows a 3D slice in $(\pos_x,\pos_y,\phi)$ space at $\vel_x=\vel_y=1, \omega=0$. The yellow set represents the target set $\targetset$, with the BRS in other colors. Shown on the right in Fig. \ref{fig:Quad6D_Combo} are 3D slices in $(\vel_x,\vel_y,\omega)$ space at $\pos_x,\pos_y=1.5, \phi=1.5$ through different points in time. The sets grow darker as time propagates backward. The union of the BRSs is the BRT, shown as the gray surface.

\subsection{The 10D Near-Hover Quadrotor \label{sec:hd10}}
The 10D Near-Hover Quadrotor was used for experiments involving learning-based MPC \cite{Bouffard12}. Its dynamics are

\begin{equation}
\label{eq:Quad10D_dyn}
\begin{aligned}
\begin{array}{c}
	\left[
	\begin{array}{c}
	\dot\pos_x\\
	\dot\vel_x\\
	\dot\theta_x\\
	\dot\omega_x\\
	\dot\pos_y\\
	\dot\vel_y\\
	\dot\theta_y\\
	\dot\omega_y\\
	\dot\pos_z\\
	\dot\vel_z
	\end{array}
	\right]
	=
	\left[
	\begin{array}{c}
	\vel_x + \dstb_x\\
	g \tan \theta_x\\
	-d_1 \theta_x + \omega_x\\
	-d_0 \theta_x + n_0 S_x\\
	\vel_y + \dstb_y\\
	g \tan \theta_y\\
	-d_1 \theta_y + \omega_y\\
	-d_0 \theta_y + n_0 S_y\\
	\vel_z + \dstb_z \\
	k_T T_z - g
	\end{array}
	\right]
\end{array}\\
\end{aligned}
\end{equation}

\noindent where $(\pos_x, \pos_y, \pos_z)$ denotes the position, $(\vel_x, \vel_y, \vel_z)$ denotes the velocity, $(\theta_x, \theta_y)$ denotes the pitch and roll, and $(\omega_x, \omega_y)$ denotes the pitch and roll rates. The controls of the system are $(S_x, S_y)$, which respectively represent the desired pitch and roll angle, and $T_z$, which represents the vertical thrust. The system experiences the disturbance $(\dstb_x, \dstb_y, \dstb_z)$ which represents wind in the three axes. $g$ denotes the acceleration due to gravity. The parameters $d_0, d_1, n_0, k_T$, as well as the control bounds $\cset$, that we used were $d_0 = 10, d_1 = 8, n_0 = 10, k_T = 0.91, |\ctrl_x|, |\ctrl_y| \le 10 \text{ degrees}, 0 \le \ctrl_z \le 2g, |\dstb_x|, \dstb_y \le 0.5\text{ m/s}, |\dstb_z| \le 1\text{ m/s}$. The system can be fully decoupled into three subsystems of 4D, 4D, and 2D, respectively:

\begin{equation}
\scstate_1 = (\pos_x, \vel_x, \theta_x, \omega_x), \scstate_2 = (\pos_y, \vel_y, \theta_y, \omega_y), \scstate_3 = (\pos_z, \vel_z)
\end{equation}

The target set is chosen to be

\begin{equation}
\begin{aligned}
\targetset = \{(&\pos_x, \vel_x, \theta_x, \omega_x, \pos_y, \vel_y, \theta_y, \omega_y, \pos_z, \vel_z): \\
&\qquad |\pos_x|, |\pos_y| \le 1, |\pos_z| \le 2.5\}
\end{aligned}
\end{equation}

This target set can be written as $\targetset = \bigcap_{i=1}^3 \bp(\targetset_{\scstate_i})$, where $\bp(\targetset_{\scstate_i}), i = 1, 2, 3$ are given by

\begin{equation}
\begin{aligned}
\sctarget{1} &= \{(\pos_x, \vel_x, \theta_x, \omega_x): |\pos_x| \le 1\} \\
\sctarget{2} &= \{(\pos_y, \vel_y, \theta_y, \omega_y): |\pos_y| \le 1\} \\
\sctarget{3} &= \{(\pos_z, \vel_z): |\pos_z| \le 2.5\}
\end{aligned}
\end{equation}

Since the subsystems do not have any common controls or disturbances, and $\targetset = \bigcap_{i=1}^3 \bp(\targetset_{\scstate_i})$, we can compute the full-dimensional $\maxbrs(t)$ and $\maxbrt(t)$ by reconstructing lower-dimensional BRSs and BRTs. A discussion of disturbances can be found in Section \ref{sec:dstb}.

\begin{figure}[H]
  \centering
  \includegraphics[width=\columnwidth]{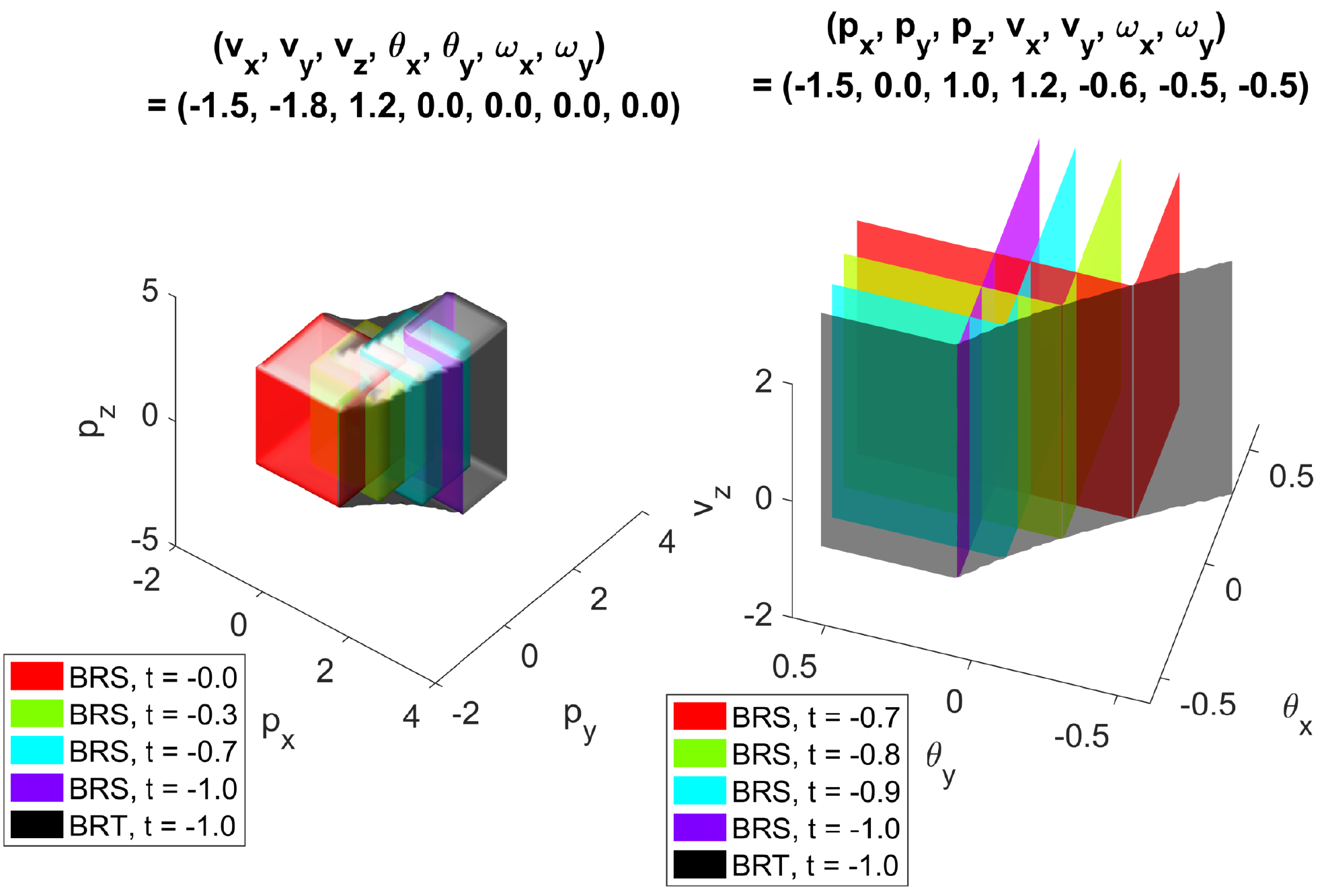}
  \caption{3D slices of the 10D BRSs over time (colored surfaces) and BRT (black surface) for the Near-Hover Quadrotor. The slices are taken at the indicated 7D point.}
  \label{fig:Quad10D}	
\end{figure}

From the target set, we computed the 10D BRS and BRT, $\maxbrs(s), \maxbrt(s), s\in [-1, 0]$. In the left subplot of Fig. \ref{fig:Quad10D}, we show a 3D slice of the BRS and BRT sliced at $(\vel_x, \vel_y, \vel_z) = (-1.5, -1.8, 1.2), \theta_x = \theta_y = \omega_x = \omega_y = 0$. The colored sets show the slice of the BRSs $\maxbrs(s), s \in [-1, 0]$, with the times color-coded according to the legend. The slice of the BRT is shown as the black surface; the BRT is the union of BRSs by Proposition \ref{prop:sets2tube_goal}.


The right subplot of Fig. \ref{fig:Quad10D} shows the BRS and BRT in $(\theta_x, \theta_y, \vel_z)$ space, sliced at $(\pos_x, \pos_y, \pos_z) = (-1.5, 0, 1), (\vel_x, \vel_y) = (1.2, -0.6), \omega_x = \omega_y = -0.5$. To the best of our knowledge, such a slice of the exact BRS and BRT is not possible to obtain using previous methods, since a high-dimensional system model like \eqref{eq:Quad10D_dyn} is needed for analyzing the angular behavior of the system.  

\section{Handling disturbances\label{sec:dstb}}
Under the presence of disturbances, the full system dynamics changes from \eqref{eq:fdyn} to

\begin{equation}
\begin{aligned}
\label{eq:fdyn_dstb}
\frac{d\state}{ds} = \dot\state = \fdyn(\state, \ctrl, \dstb), s \in [t, 0], \ctrl \in \cset, \dstb \in \dset
\end{aligned}
\end{equation}

\noindent where $\dstb \in \dset$ represents the disturbance, with $\dstb(\cdot) \in \dfset$ drawn from the set of measurable functions.
%

In addition, we assume that the disturbance function $\dstb(\cdot)$ is drawn from the set of non-anticipative strategies \cite{Mitchell05}, denoted $\Gamma(t)$. We denote the mapping from $\ctrl(\cdot)$ to $\dstb(\cdot)$ as $\gamma[\ctrl](\cdot)$ as in \cite{Mitchell05}. The subsystems in \eqref{eq:scdyn} are now written as


\begin{equation}
\label{eq:scdyn_dstb}
\begin{aligned}
\dot\spart_1 &= \fdyn_1(\spart_1, \spart_c, \ctrl, \dstb) \\
\dot\spart_2 &= \fdyn_2(\spart_2, \spart_c, \ctrl, \dstb) \\
\dot\spart_c &= \fdyn_c(\spart_c, \ctrl, \dstb) \\
\end{aligned}
\end{equation}

In general, subsystem disturbances may have shared components. Whether this is the case is very important, as some of the results involving disturbances become stronger when the subsystem disturbances do not have shared components.

Trajectories of the system and subsystems are now denoted $\traj(s; \state, t, \ctrl(\cdot), \dstb(\cdot)), \sctraj_i(s; \scstate_i, t, \ctrl(\cdot), \dstb(\cdot))$, and satisfy conditions analogous to \eqref{eq:fdyn_traj} and \eqref{eq:scdyn_traj} respectively. We also need to incorporate the disturbance into the BRS and BRT definitions:

\begin{equation}
\begin{aligned}
\minbrt(t) &= \{\state: \exists\gamma[\ctrl](\cdot), \forall \ctrl(\cdot), \exists s \in [t, 0], \\
&\qquad\qquad\qquad \traj(s; \state, t, \ctrl(\cdot), \gamma[\ctrl](\cdot)) \in \targetset \} \\
\maxbrt(t) &= \{\state: \forall\gamma[\ctrl](\cdot),  \exists \ctrl(\cdot), \exists s \in [t, 0], \\
&\qquad\qquad\qquad \traj(s; \state, t, \gamma[\ctrl](\cdot), \dstb(\cdot)) \in \targetset \} \\
\minbrs(t) &= \{\state: \exists\gamma[\ctrl](\cdot), \forall \ctrl(\cdot), \traj(0; \state, t, \ctrl(\cdot), \gamma[\ctrl](\cdot)) \in \targetset \} \\
\maxbrs(t) &= \{\state: \forall\gamma[\ctrl](\cdot), \exists \ctrl(\cdot),  \traj(0; \state, t, \ctrl(\cdot), \gamma[\ctrl](\cdot)) \in \targetset \}
\end{aligned}
\end{equation}

Subsystem BRSs $\maxbrs_i, \minbrs_i, i = 1,2$ are defined analogously.

\subsection{Self-Contained Subsystems}
\label{sec:sc_dstb}
Under the presence of disturbances, the results from Section \ref{sec:sc} carry over with some modifications. Theorems \ref{thm:sc_reach_un} and \ref{thm:sc_reach_in} need to be changed slightly, and the reconstructed BRS is now an approximation conservative in the right direction.

\begin{cor}	\label{cor:maxBRS_Union_SCS_Disturbance}
	Suppose that the full system in \eqref{eq:fdyn} can be decomposed into the form of \eqref{eq:scdyn_dstb}, then

\begin{equation}
\label{eq:maxBRS_union_disturbance}
	\begin{aligned}
&\targetset = \bp(\sctarget{1}) \cup \bp(\sctarget{2}) \\
&\Rightarrow \maxbrs(t) \supseteq \bp(\maxbrs_1(t)) \cup \bp(\maxbrs_2(t))
\end{aligned}
\end{equation}

To solve this, in the proof of Theorem \ref{thm:sc_reach_un}, \eqref{eq:subsystem_reach_target} becomes

\begin{equation}
\label{eq:dstb_statement}
\begin{aligned}
\forall\gamma[\ctrl](\cdot), \exists \ctrl(\cdot), &\quad\sctraj_1(0; \bar\scstate_1, t, \ctrl(\cdot), \dstb(\cdot)) \in \sctarget{1} \vee\\
&\quad\sctraj_2(0; \bar\scstate_2, t, \ctrl(\cdot), \dstb(\cdot)) \in \sctarget{2}
\end{aligned}
\end{equation}

The expression ``$\forall\gamma[\ctrl](\cdot), \exists \ctrl(\cdot)$'' can no longer be distributed, thus making the reconstructed BRS a conservative approximation of the true BRS in the right direction. By conservative in the right direction, we mean that a state $\state$ in the reconstructed BRS is guaranteed to be able to reach the target.
\end{cor}
%
%
%
%
%
%

\begin{cor} \label{cor:minBRS_Intersection_SCS_Disturbance}
	  Suppose that the full system in \eqref{eq:fdyn} can be decomposed into the form of \eqref{eq:scdyn_dstb}, then
	
	\begin{equation}
	\label{eq:minBRS_intersection_disturbance}
	\begin{aligned}
	&\targetset = \bp(\sctarget{1}) \cap \bp(\sctarget{2}) \\
	&\Rightarrow \minbrs(t) \subseteq \bp(\minbrs_1(t)) \cap \bp(\minbrs_2(t))
	\end{aligned}
	\end{equation}
	
	The proof of Theorem \ref{thm:sc_reach_in} makes the same arguments except that it involves complements of sets instead. Again, the reconstructed BRS is a conservative approximation of the true BRS in the right direction, meaning that a state $\state$ outside of the reconstructed BRS is guaranteed to be able to avoid the target.
\end{cor}

If the subsystem disturbances have no shared components, then \eqref{eq:dstb_statement} becomes

\begin{equation*}
\begin{aligned}
\forall \left(\gamma_1[\ctrl](\cdot), \gamma_2[\ctrl](\cdot)\right),\exists \ctrl(\cdot), &\quad\sctraj_1(0; \bar\scstate_1, t, \ctrl(\cdot), \gamma_1[\ctrl](\cdot)) \in \sctarget{1} \vee\\
&\quad\sctraj_2(0; \bar\scstate_2, t, \ctrl(\cdot), \gamma_2[\ctrl](\cdot)) \in \sctarget{2}
\end{aligned}
\end{equation*}

\noindent where $\gamma[\ctrl](\cdot)$ is written as $\left(\gamma_1[\ctrl](\cdot), \gamma_2[\ctrl](\cdot)\right)$. 

In this case, the expression ``$\forall \left(\gamma_1[\ctrl](\cdot), \gamma_2[\ctrl](\cdot)\right),\exists \ctrl(\cdot)$'' can be distributed. Therefore, in this case Theorems \ref{thm:sc_reach_un} and \ref{thm:sc_reach_in} still hold.

\subsection{Subsystems with Decoupled Control}
\label{sec:decoupled_dstb}
For systems with decoupled control, but shared disturbance in the subsystems, results from Section \ref{sec:sc_dstb} still hold since the system dynamics structure is a special case of that in Section \ref{sec:sc_dstb}. In addition, results from Section \ref{sec:decoupled} hold with some modifications. Propositions \ref{prop:decoupled_reach_in} and \ref{prop:decoupled_reach_un} need to be modified, and again the reconstructed BRS is now an approximation conservative in the right direction. 

\begin{cor}
	\label{cor:targetset_in_decoupled_with_dstb}
	Suppose that the full system in \eqref{eq:fdyn} can be decomposed into the form of \eqref{eq:ddyn}, with the addition of shared disturbances. Then,
\begin{equation}
\label{eq:targetset_in_decoupled_with_dstb}
\begin{aligned}
&\targetset = \bp(\sctarget{1}) \cap \bp(\sctarget{2}) \\
&\Rightarrow \maxbrs(t) \supseteq \bp(\maxbrs_1(t)) \cap \bp(\maxbrs_2(t))
\end{aligned}
\end{equation}
To prove this, we modify Proposition \ref{prop:decoupled_reach_in} by changing \eqref{eq:decoupled_ctrl_before_distr} to
\begin{equation}
\begin{aligned}
&\forall \gamma[\ctrl](\cdot), \exists \left(\cpart_1(\cdot), \cpart_2(\cdot)\right) \\
&\quad\big(\eta_1(s; \bar\spart_1, \bar\spart_c, t, \cpart_1(\cdot), \gamma[\ctrl](\cdot)), \eta_c(s; \bar\spart_c, t)\big) \in \sctarget{1} \wedge\\
&\quad\big(\eta_2(s; \bar\spart_2, \bar\spart_c, t, \cpart_2(\cdot), \gamma[\ctrl](\cdot)), \eta_c(s; \bar\spart_c, t)\big) \in \sctarget{2}
\end{aligned} \\
\end{equation}
The expression ``$\forall \gamma[\ctrl](\cdot), \exists \left(\cpart_1(\cdot), \cpart_2(\cdot)\right)$'' cannot be distributed to lead to a statement analogous to \eqref{eq:decoupled_ctrl_after_distr}. Hence, the forward direction of Proposition \ref{prop:decoupled_reach_in} does not hold, and conservativeness is introduced.
\end{cor}

By the same reasoning, the result of Proposition \ref{prop:decoupled_reach_un} changes to the following.

\begin{cor}
	\label{cor:targetset_un_decoupled_with_dstb}
\begin{equation}
\label{eq:targetset_un_decoupled_with_dstb}
\begin{aligned}
&\targetset = \bp(\sctarget{1}) \cup \bp(\sctarget{2}) \\
&\Rightarrow \minbrs(t) \subseteq \bp(\minbrs_1(t)) \cup \bp(\minbrs_2(t))
\end{aligned}
\end{equation}
\end{cor}

In both cases, conservative approximations of the BRS can still be obtained. 

\begin{table}[h]
	\centering
	\caption{BRS Results from Subsections \ref{sec:sc_dstb} \& \ref{sec:decoupled_dstb}}
	\SHnote{
	{\scriptsize
	\label{table:BRS_Disturbance}
	\begin{tabular}{|l|l|l|l|l|}
		\hline
		\textbf{Shared Controls}       & \multicolumn{2}{c|}{\textbf{Yes}}                                                                                                                                                                               & \multicolumn{2}{c|}{\textbf{No}}                                                                                                                                                                                               \\ \hline
		\textbf{Shared Dstb.}    & \multicolumn{2}{c|}{\textbf{Yes}}                                                                                                                                                                               & \multicolumn{2}{c|}{\textbf{Yes}}                                                                                                                                                                                              \\ \hline
		\textbf{Target}                & \textbf{Intersection}                                                                                             & \textbf{Union}                                                                                       & \textbf{Intersection}                                                                                              & \textbf{Union}                                                                                                     \\ \hline
		\begin{tabular}[c]{@{}l@{}}Recover \\ Max. BRT?  \end{tabular}     & No                                                                                                       & Yes, consrv                                                                                & Yes, consrv                                                                                              & Yes, consrv                                                                                                        \\ \hline
		\begin{tabular}[c]{@{}l@{}}Recover \\ Min. BRT?  \end{tabular}     & Yes, consrv                                                                                             & No                                                                                          & Yes, consrv                                                                                                        & Yes, consrv                                                                                              \\ \hline
		Equation(s) & Cor \ref{cor:minBRS_Intersection_SCS_Disturbance}, (\ref{eq:minBRS_intersection_disturbance}) & Cor \ref{cor:maxBRS_Union_SCS_Disturbance},  (\ref{eq:maxBRS_union_disturbance}) & Cor \ref{cor:targetset_in_decoupled_with_dstb}, (\ref{eq:targetset_in_decoupled_with_dstb}) & Cor \ref{cor:targetset_un_decoupled_with_dstb}, (\ref{eq:targetset_un_decoupled_with_dstb}) \\ \hline
	\end{tabular}
}}
\end{table}

\subsection{Decomposition of Reachable Tubes}
\label{sec:set2tube_dstb}
Under disturbances, the results from Section \ref{sec:set_to_tube} carry over with modifications. For reconstruction from other BRTs, the arguments in Proposition \ref{prop:BRT_Union} do not change. However, in the case where there are overlapping components in the subsystem disturbances, the reconstructed BRTs become conservative approximations:

\begin{cor} Suppose our system has coupled control and disturbance as in \eqref{eq:scdyn_dstb}, then
	\label{cor:maxBRT_union_disturbance}
	\begin{equation}
	\label{eq:maxBRT_union_disturbance}
	  \begin{aligned}
	&\targetset = \bp(\sctarget{1}) \cup \bp(\sctarget{2}) \\
	&\Rightarrow \maxbrt(t) \supseteq \bp(\maxbrt_1(t)) \cup \bp(\maxbrt_2(t))
	\end{aligned}
\end{equation}
\end{cor}

\begin{cor} Suppose our system has subsystem controls that do not share any components, then
	\label{cor:minBRT_union_disturbance}
		\begin{equation}
	\label{eq:minBRT_union_disturbance}
	\begin{aligned}
	&\targetset = \bp(\sctarget{1}) \cup \bp(\sctarget{2}) \\
	&\Rightarrow \minbrt(t) \subseteq \bp(\minbrt_1(t)) \cup \bp(\minbrt_2(t))
	\end{aligned}
\end{equation}
\end{cor}

For Proposition \ref{prop:sets2tube_goal}, the union of the BRSs now becomes an under-approximation of the BRT in general:
\begin{cor}Suppose our system has coupled control and disturbance as in \eqref{eq:scdyn_dstb}, then
	\label{cor:maxBRT_intersection_disturbance}
	\begin{equation}
	\label{eq:maxBRT_intersection_disturbance}
	\begin{aligned}
	&\targetset = \bp(\sctarget{1}) \cap \bp(\sctarget{2}) \\
	&\Rightarrow \bigcup_{s \in [t, 0]} \maxbrs(s) \subseteq \maxbrt(t)
	\end{aligned}
	\end{equation}
	
	 To show this, all arguments in the proof of Proposition \ref{prop:sets2tube_goal} remain the same, except \eqref{eq:union_goal_expression} no longer implies \eqref{eq:goal_after_swap}. Instead, the implication is unidirectional:
	
	\begin{equation}
	\begin{aligned}
	&\exists s \in [t, 0], \forall \gamma[\ctrl](\cdot), \exists \ctrl(\cdot), \traj(0; \state, s, \ctrl(\cdot), \gamma[\ctrl](\cdot)) \in \targetset \\
	&\Rightarrow \forall \gamma[\ctrl](\cdot), \exists \ctrl(\cdot), \exists s \in [t, 0], \traj(s; \state, t, \ctrl(\cdot), \gamma[\ctrl](\cdot)) \in \targetset
	\end{aligned}
	\end{equation}
	
	This is due to the switching of the order of the expressions ``$\exists s \in [t, 0]$'' and ``$\gamma[\ctrl](\cdot)$''. Therefore, the union of the BRSs becomes an under-approximation of the BRT, a conservatism in the right direction: a state in the under-approximated BRT is still guaranteed to be able to reach the target.
\end{cor}

In contrast to Proposition \ref{prop:sets2tube_goal}, all the arguments of Theorem \ref{thm:sets2tube_avoid} hold, since there no change of order of expressions involving existential and universal quantifiers.

\begin{table}[h]
	\centering
	\caption{BRT Results for Reconstruction from Tubes}
	\label{table:BRT_Tubes_Dstb}
	\SHnote{{\scriptsize
		\begin{tabular}{|l|c|l|l|l|}
			\hline
			\textbf{Shared Controls}    & \multicolumn{2}{c|}{\textbf{Yes}}                         & \multicolumn{2}{c|}{\textbf{No}}     \\ \hline
			\textbf{Shared Dstb.} & \multicolumn{2}{c|}{\textbf{Yes}}                         & \multicolumn{2}{c|}{\textbf{Yes}}    \\ \hline
			\textbf{Target}             & \multicolumn{1}{l|}{\textbf{Intersection}} & \textbf{Union}        & \textbf{Intersection} & \textbf{Union}        \\ \hline
			\begin{tabular}[c]{@{}l@{}}Recover \\ Max. BRT?  \end{tabular}  & No                                & Yes, consrv & No           & Yes, consrv \\ \hline
			\begin{tabular}[c]{@{}l@{}}Recover \\ Min. BRT?  \end{tabular} & No                                & No & No           & Yes, consrv \\ \hline
			Equation(s)        & N/A                               & Cor \ref{cor:maxBRT_union_disturbance}, (\ref{eq:maxBRT_union_disturbance})            & N/A          & \begin{tabular}[c]{@{}l@{}}Cor \ref{cor:maxBRT_union_disturbance}, (\ref{eq:maxBRT_union_disturbance})\\ Cor \ref{cor:minBRT_union_disturbance}, (\ref{eq:minBRT_union_disturbance})\end{tabular}            \\ \hline
	\end{tabular}}}
\end{table}

\begin{table}[h]
	\centering
	\caption{BRT Results for Reconstruction from Sets}
\label{table:BRT_Sets_Dstb}
\SHnote{{\scriptsize
	\begin{tabular}{|l|c|l|c|l|}
		\hline
		\textbf{Shared Controls}    & \multicolumn{2}{c|}{\textbf{Yes}}                  & \multicolumn{2}{c|}{\textbf{No}}                   \\ \hline
		\textbf{Shared Disturbance} & \multicolumn{2}{c|}{\textbf{Yes}}                   & \multicolumn{2}{c|}{\textbf{Yes}}                   \\ \hline
		\textbf{Target}             & \multicolumn{1}{l|}{\textbf{Intersection}} & \textbf{Union} & \multicolumn{1}{l|}{\textbf{Intersection}} & \textbf{Union} \\ \hline
		Recover Max. BRT?  & \multicolumn{4}{c|}{Yes, conserv}                                                     \\ \hline
		Recover Min. BRT?  & \multicolumn{4}{c|}{Yes, exact*}                                                      \\ \hline
		Equation(s)        & \multicolumn{4}{c|}{\begin{tabular}[c]{@{}l@{}}Cor \ref{cor:maxBRT_intersection_disturbance}, (\ref{eq:maxBRT_intersection_disturbance})\\ Thm \ref{thm:sets2tube_avoid}, (\ref{eq:minBRT_intersection})\end{tabular}}                                                           \\ \hline
	\end{tabular}
}}
	\begin{flushleft}
		* the solution here can be found only if the minimum BRSs are non-empty for the entire time period.
	\end{flushleft}
\end{table}

\subsection{Dubins Car with Disturbances}
Under disturbances, the Dubins Car dynamics are given by

\begin{equation}
\begin{aligned}
\left[ \begin{array}{c}
\dot\pos_x\\
\dot\pos_y \\
\dot\theta
\end{array} \right]
=
\left[
\begin{array}{c}
v \cos\theta + \dstb_x\\
v \sin\theta + \dstb_y\\
\omega + \dstb_\theta
\end{array}\right] \\
\omega \in \cset, \quad (\dstb_x, \dstb_y, \dstb_\theta) \in \dset
\end{aligned}
\end{equation}

\noindent with state $\state = (\pos_x, \pos_y, \theta)$, control $\ctrl = \omega$, and disturbances $\dstb = (\dstb_x, \dstb_y, \dstb_\theta)$. The state partitions are $\spart_1 = \pos_x, \spart_2 = \pos_y, \spart_c = \theta$. The subsystems states $\scstate_i$, controls $\scctrl_i$, and disturbances $\scdstb_i$ are
\begin{equation}
\begin{aligned}
\dot{\scstate_1} = 
\left[ \begin{array}{c}
\dot{\spart_1}\\
\dot{\spart_c}
\end{array} \right]
=
\left[ \begin{array}{c}
\dot{\pos_x}\\
\dot{\theta}
\end{array} \right]
&=
\left[\begin{array}{c}
v \cos\theta + \dstb_x\\
\omega +\dstb_\theta\\
\end{array}\right]\\
\dot{\scstate_2} = 
\left[ \begin{array}{c}
\dot{\spart_2}\\
\dot{\spart_c}
\end{array} \right]
=
\left[ \begin{array}{c}
\dot{\pos_y}\\
\dot{\theta}
\end{array} \right]
&=
\left[\begin{array}{c}
v \sin\theta +\dstb_y \\
\omega + \dstb_\theta
\end{array}\right]\\
\cpart_c &= \omega = \ctrl \\
\dpart_1 &= \dstb_x, \dpart_2 = \dstb_y, \dpart_c = \dstb_\theta
\end{aligned}
\end{equation}

\noindent where the overlapping state is $\theta = \spart_c$. We assume that each component of disturbance is bounded in some interval centered at zero: $|\dstb_x| \le \bar\dstb_x,|\dstb_y| \le \bar\dstb_y,|\dstb_\theta| \le \bar\dstb_\theta$. The subsystem disturbances $\scdstb_1$ and $\scdstb_2$ have the shared component $\dstb_\theta$.

\begin{figure}[H]
	\centering
	\includegraphics[width=0.9\columnwidth]{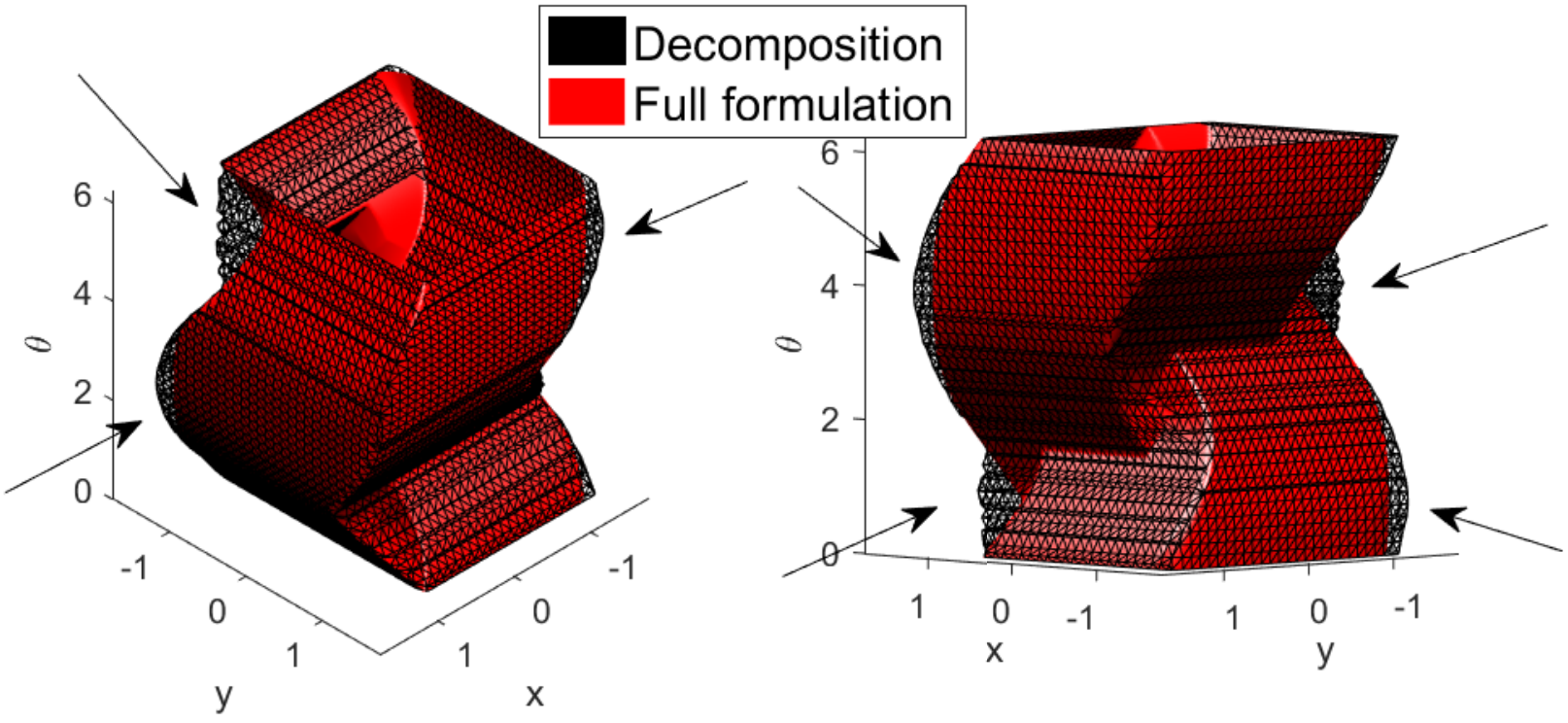}
	\caption{Minimal BRTs computed directly in 3D and via decomposition in 2D for the Dubins Car under disturbances with shared components. The reconstructed BRT is an over-approximation of the true BRT. The over-approximated regions of the reconstruction are indicated by the arrows.}
	\label{fig:dubins_dstb_coupled}
\end{figure}

\begin{figure}[H]
  \centering
  \includegraphics[width=0.9\columnwidth]{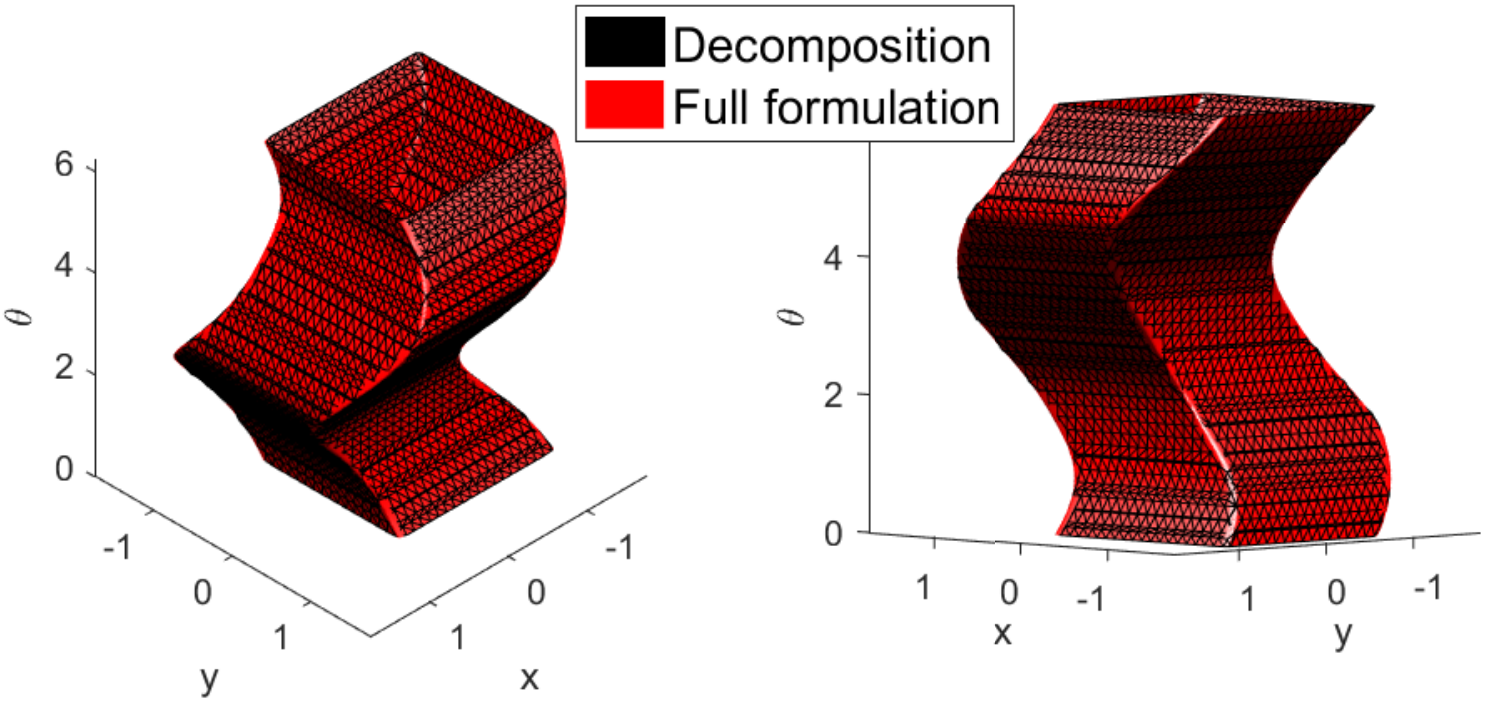}
  \caption{Minimal BRTs computed directly in 3D and via decomposition in 2D for the Dubins Car under disturbances \textit{without} shared components. In this case, the BRT computed using decomposition matches the true BRT.}
  \label{fig:dubins_dstb_uncoupled}
\end{figure}

Fig. \ref{fig:dubins_dstb_coupled} compares the BRT $\minbrt(t), t = -0.5$ computed directly from the target set in \eqref{eq:dubins_target}, and using our decomposition technique from the subsystem target sets in \eqref{eq:dubins_target_decomp}. For this computation, we chose $\bar\dstb_x, \bar\dstb_y = 1, \bar\dstb_\theta = 5$. 

Since there is a shared component in the disturbances, the BRT computed using our decomposition technique becomes an over-approximation of the true BRT. One can see the over-approximation by noting that the black set is not flush against the red set, as marked by the arrows in Fig. \ref{fig:dubins_dstb_coupled}.

Fig. \ref{fig:dubins_dstb_uncoupled} shows the same computation with $\bar\dstb_\theta = 0$, so that subsystem disturbances effectively have no shared components. In this case, one can see that the BRTs computed directly in 3D and via decomposition in 2D are the same.


\section{Conclusions and Future Work}
In this paper, we presented a general system decomposition method for efficiently computing BRSs and BRTs in several scenarios. By performing computations in lower-dimensional subspaces, computation burden is substantially reduced, allowing currently tractable computations to be orders of magnitude more faster, and currently intractable computations to become tractable. Unlike related work on computation of BRSs and BRTs, our method can significantly reduce dimensionality without sacrificing any optimality. 

Under disturbances, the reconstructed BRSs and BRTs sometimes become slightly conservative approximations which are still useful for providing performance and safety guarantees. To the best of our knowledge, such guarantees for high-dimensional systems are now possible for the first time. Our decomposition technique can also be used in combination with other dimensionality reduction or approximation techniques, further alleviating the curse of dimensionality.

We are currently extending our decomposition technique to other scenarios, including more representations of full-dimensional sets in lower-dimensional subspaces and more families of system dynamics. In addition, we look forward to combining our technique with other related techniques such as reinforcement learning and machine learning, automating the system decomposition process, and demonstrating our theory in hardware experiments.



\bibliographystyle{IEEEtran}
\bibliography{IEEEabrv,references}

\begin{IEEEbiography}[{\includegraphics[width=1in,height=1.25in,clip,keepaspectratio]{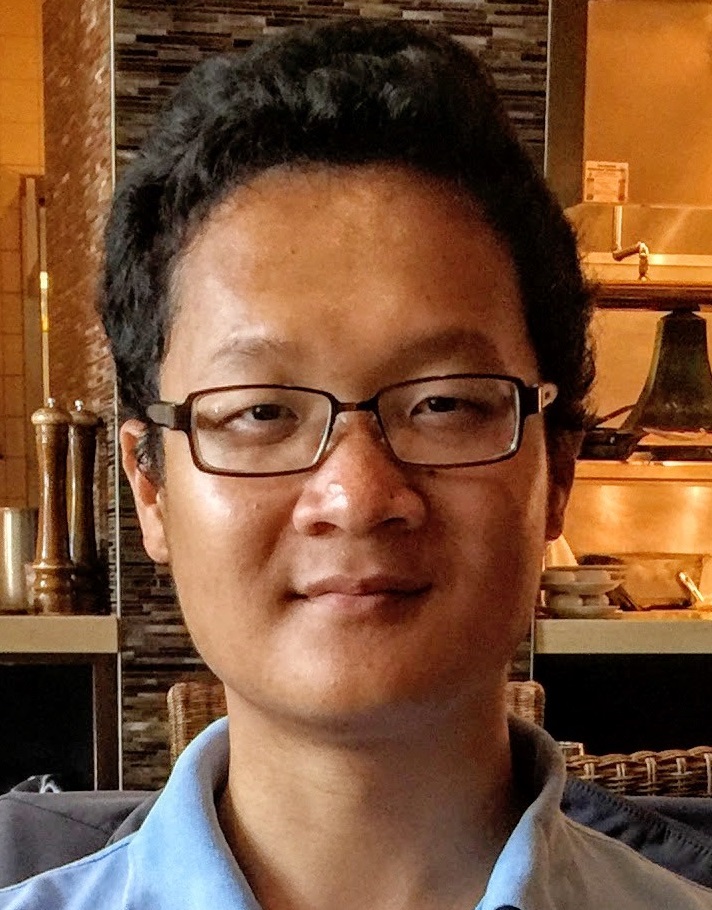}}]{Mo Chen}
received the B.A.Sc. degree in Engineering Physics from the University of British Columbia, Vancouver, BC, Canada, in 2011, and is currently a Ph.D. candidate in Electrical Engineering and Computer Sciences at the University of California, Berkeley.
\end{IEEEbiography}
\begin{IEEEbiography}[{\includegraphics[width=1in,height=1.25in,clip,keepaspectratio]{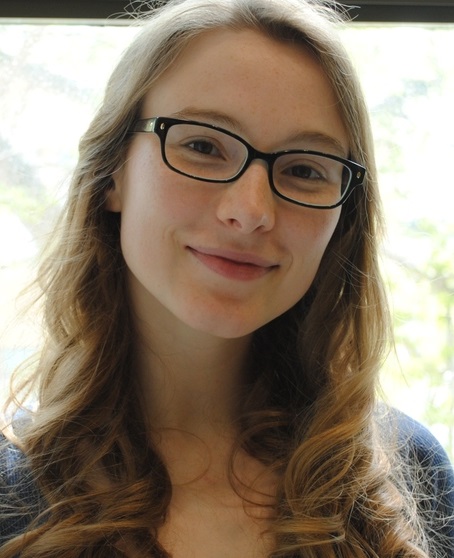}}]{Sylvia L. Herbert}
received her B.S. and M.S. degrees in Mechanical Engineering at Drexel University, Philadelphia, PA, in 2014. She is currently a Ph.D. student in Electrical Engineering and Computer Sciences at the University of California, Berkeley.
\end{IEEEbiography}
\begin{IEEEbiography}[{\includegraphics[width=1in,height=1.25in,clip,keepaspectratio]{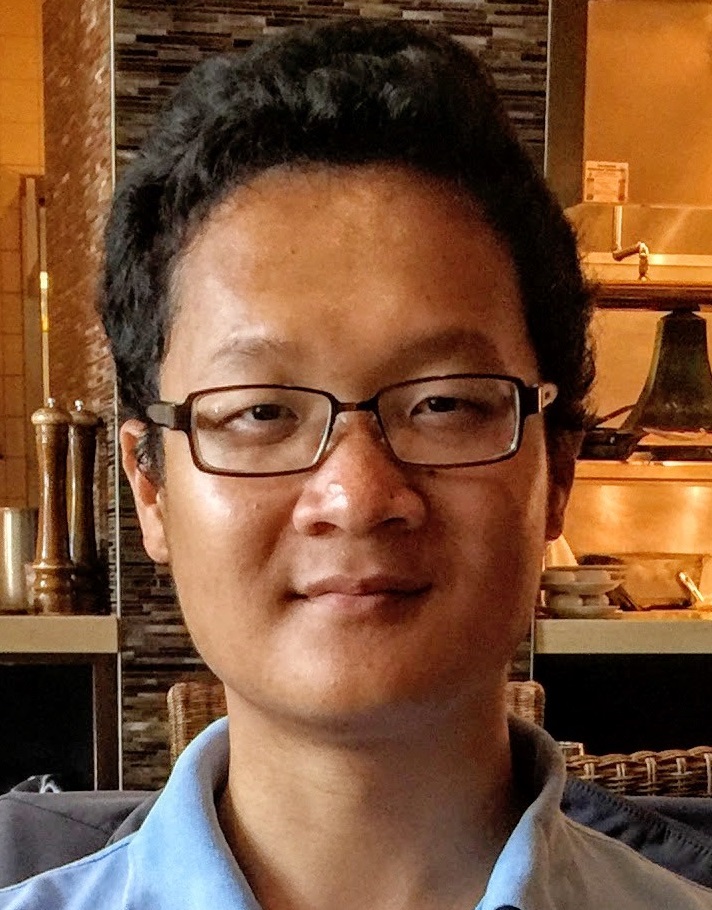}}]{Mahesh S. Vashishtha}
is currently an undergraduate student in Computer Science at the University of California, Berkeley
\end{IEEEbiography}
\begin{IEEEbiography}[{\includegraphics[width=1in,height=1.25in,clip,keepaspectratio]{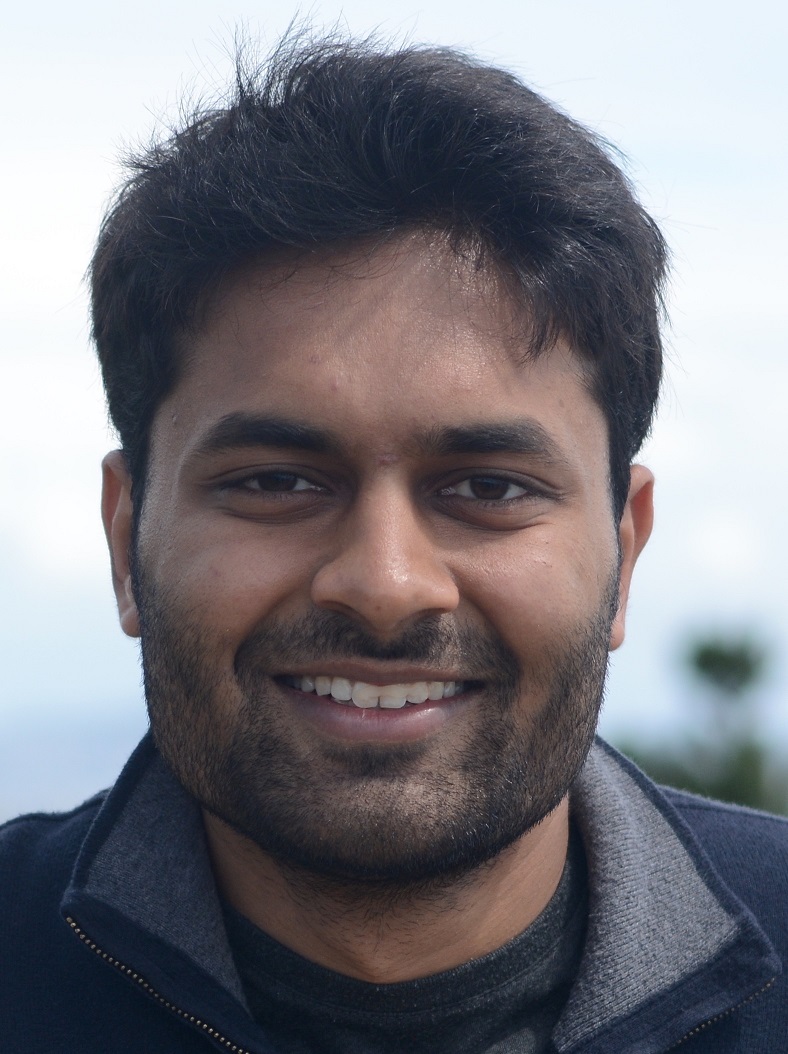}}]{Somil Bansal}
received the B.Tech. degree in Electrical Engineering from the Indian Institute of Technology, Kanpur, India, in 2012, and is currently a Ph.D. student in Electrical Engineering and Computer Sciences at the University of California, Berkeley.
\end{IEEEbiography}
\begin{IEEEbiography}[{\includegraphics[width=1in,height=1.25in,clip,keepaspectratio]{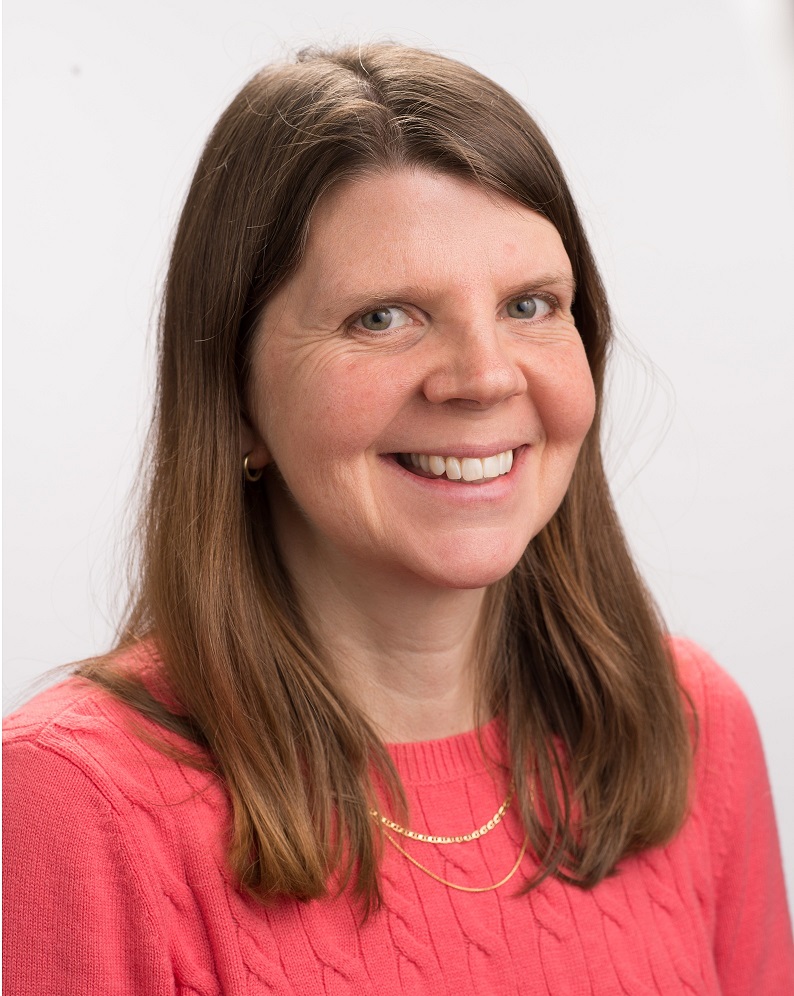}}]{Claire J. Tomlin}
is the Charles A. Desoer Professor of Engineering in Electrical Engineering and Computer Sciences at the University of California, Berkeley. She was an Assistant, Associate, and Full Professor in Aeronautics and Astronautics at Stanford from 1998 to 2007, and in 2005 joined Berkeley. Claire works in the area of control theory and hybrid systems, with applications to air traffic management, UAV systems, energy, robotics, and systems biology.  She is a MacArthur Foundation Fellow (2006) and an IEEE Fellow (2010), and in 2010 held the Tage Erlander Professorship of the Swedish Research Council at KTH in Stockholm.
\end{IEEEbiography}
\end{document}